\documentclass{amsart}

\usepackage{amsmath,amssymb,amsthm,multicol,mathtools,hyperref,dsfont,pinlabel,enumitem,dsfont}
\usepackage{comment}
\usepackage[usenames,dvipsnames]{xcolor}
\usepackage{tikz-cd}
\usepackage{mathrsfs}
\usepackage[numbers]{natbib}

\newcommand\bs{\boldsymbol}
\newcommand\cut{\setminus\!\setminus}
\newcommand\wt[1]{\scalebox{.9}{$\widetilde{
#1
}$}}

\newcommand\x{\textbf{x}}

\newcommand\Z{\mathds{Z}}

\newcommand\ess{\text{ess}}

\newcommand\Caps{\text{Caps}}

\newcommand\bbm{\begin{bmatrix}}
\newcommand\ebm{\end{bmatrix}}

\newcommand\White[1]{\color{white}#1\color{black}}

\definecolor{myPurple}{cmyk}{.5,.75,0,0}
\definecolor{myBlue}{cmyk}{.5,1,0,0}
\definecolor{myGreen}{cmyk}{1,0,1,0}
\definecolor{myGray}{cmyk}{0,0,0,.5}
\definecolor{myPink}{cmyk}{0,.73,.2,0}
\definecolor{myRed}{cmyk}{0,1,1,0}

\newcommand\MobPos{\raisebox{-2pt}{\includegraphics[height=11pt]{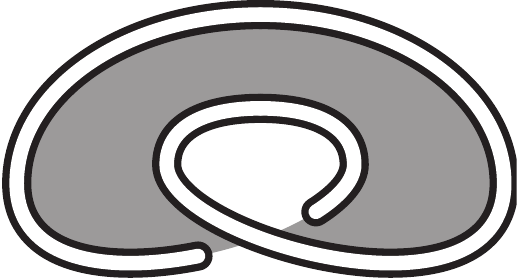}}}
\newcommand\MobNeg{\raisebox{-2pt}{\includegraphics[height=11pt]{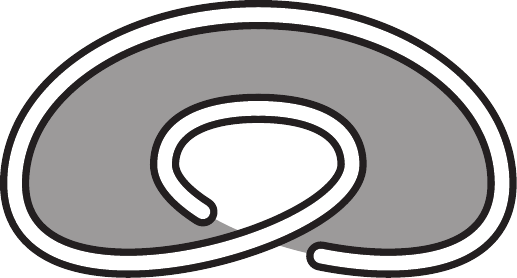}}}

\theoremstyle{plain}
\newtheorem{theorem}{Theorem}[section]

\newtheorem{obs}[theorem]{Observation}
\newtheorem{prop}[theorem]{Proposition}
\newtheorem{cor}[theorem]{Corollary}

\newtheorem*{P:GeomEss}{Proposition \ref{P:GeomEss}}
\newtheorem*{T:BadPlumb}{Theorem \ref{T:BadPlumb}}
\newtheorem*{T:PlumbEss}{Theorem \ref{T:PlumbEss}}
\newtheorem*{T:PlumbEndEss}{Theorem \ref{T:PlumbEndEss}}
\newtheorem*{T:Endess}{Theorem \ref{T:Endess}}
\newtheorem*{T:CBEss}{Theorem \ref{T:CBEss}}
\newtheorem*{C:StateEss}{Corollary \ref{C:StateEss}}
\newtheorem*{T:Singular}{Theorem \ref{T:Singular}}
\newtheorem*{T:CC}{Theorem \ref{T:CC}}
\newtheorem*{T:canon}{Theorem \ref{T:canon}}
\newtheorem*{T:DeplumbUnique}{Theorem \ref{T:DeplumbUnique}}
\newtheorem*{C:DeplumbUnique1}{Corollary \ref{C:DeplumbUnique1}}
\newtheorem*{C:DeplumbUnique2}{Corollary \ref{C:DeplumbUnique2}}
\newtheorem*{C:EssGraph}{Corollary \ref{C:EssGraph}}

\theoremstyle{definition}
\newtheorem{convention}[theorem]{Convention}
\newtheorem{notation}[theorem]{Notation}
\newtheorem{definition}[theorem]{Definition}
\newtheorem{question}[theorem]{Question}

\newtheorem{rem}[theorem]{Remark}

\numberwithin{equation}{section}

\title[How essential is a spanning surface?]{How essential is a spanning surface?}

\author{Thomas Kindred}

\address{
Department of Mathematical Sciences \\ 
Burton Hall 115, Smith College \\
Northampton, MA 01063}
\email{tkindred@smith.edu}

\begin{document}

\begin{abstract}
We introduce new numerical invariants of a spanning surface $F\subset S^3$, called the {\it algebraic essence} and {\it geometric essence} of $F$, which measure how far $F$ is from being compressible.  We extend a theorem of Ozawa by showing that algebraic essence behaves well under Murasugi sum.
We also introduce a ``twisted'' generalization of Murasugi sum and use it to compute algebraic and geometric essence for many examples, including checkerboard surfaces from reduced alternating diagrams.
Finally, we extend all of these results to Murasugi sums and twisted Murasugi sums of spanning surfaces in arbitrary 3-manifolds.
\end{abstract}



\maketitle

\section{Introduction}

A great deal of information about a knot or link $L$ in a 3-manifold $M$ is captured by the spanning surfaces for $L$.  These are embedded surfaces $F\subset M$ with $\partial F=L$.  Because every link has infinitely many distinct spanning surfaces and the simplest ones are the most likely to convey pertinent information about the link, it is common to restrict one's attention to {\it essential} surfaces.  There are, roughly, two notions of essential spanning surfaces in the literature: an ``algebraic'' notion, defined in terms of maps on fundamental groups, and a ``geometric'' notion which can be understood in terms of surgery operations. For more detail, see Definitions \ref{D:GeomEss} and \ref{D:AlgEss}, \textsection1 of \cite{ht}, or \textsection2 and the Appendix of \cite{natural}.

In this paper, we introduce the  (algebraic) {\it essence} $\ess(F)$ of a spanning surface $F$,  which measures how far a surface is from being algebraically (or $\pi_1$-) inessential 
and generalizes $\pi_1$-essentiality in the sense that $F$ is $\pi_1$-essential if and only if $\text{ess}(F)\geq 2$. 
Our main result is that this new invariant behaves well under an important operation called Murasugi sum, which involves gluing together two spanning surfaces along a disk so as to obtain another spanning surface (there is one extra condition---see Definition \ref{D:Plumb}):

\begin{T:PlumbEss}
If $F=F_0*F_1$ is a Murasugi sum of $\pi_1$-essential spanning surfaces in $S^3$, then $\ess(F)\geq \min_{i=0,1}\ess(F_i)$.
\end{T:PlumbEss}

Theorem \ref{T:PlumbEss} extends a theorem of Gabai (in the orientable case) and Ozawa \cite{gab1,gab2,ozawa11}.

We also define a related notion that we call the {\it geometric essence} of  $F$, denoted $\ess_g(F)$. A reader familiar with the notion of {\it representativity} (introduced by Pardon to study the distortion of knots \cite{pardon11}) may think of $\ess(F)$ and $\ess_g(F)$ as ways of adapting that idea to spanning surfaces. 

The main result of \cite{natural} implies Murasugi sum does not respect geometric essence, in the sense of Theorem \ref{T:PlumbEss}.
For many spanning surfaces $F$, however, algebraic and geometric essence are equal, $\ess(F)=\ess_g(F)$.  We show in particular that this is the case when $F$ is a checkerboard surface from a reduced alternating link diagram:  

\begin{T:CBEss}
The algebraic and geometric essence of any checkerboard surface from any reduced alternating diagram on $S^2\subset S^3$ both equal the length of the shortest cycle in its Tait graph.
\end{T:CBEss}

This intuitive result, however, is harder to prove than one might expect.  To this end, we introduce a ``twisted'' generalization of Murasugi sum and prove that it too respects essence in a sense analogous to Theorem \ref{T:PlumbEss}, with some extra restrictions. See Theorem \ref{T:TwistedEss} for the precise statement of that result.

Theorem \ref{T:CBEss} extends via Theorem \ref{T:PlumbEss} to a broad class of examples (see Definition \ref{D:Homogeneously_Adequate}) that includes any adequate state surface from any alternating link diagram:

\begin{C:StateEss}
Let $F_x$ be a homogeneously adequate state surface in $S^3$, and let $n$ be the length of the shortest cycle in its state graph.  Then under any layering of the state disks of $F_x$ we have $\ess_g(F_x)\geq\ess(F_x)\geq n$, and there is some layering of these disks such that $\ess_g(F_x)=\ess(F_x)=n$.
\end{C:StateEss}

We further extend these results to Murasugi sums of spanning surfaces in arbitrary 3-manifolds (also see Definitions \ref{D:EndEss} and \ref{D:Plumb3Mfld}):

\begin{T:PlumbEndEss}
Suppose $M=M_0\#M_1$ is a (possibly trivial) connect sum of 3-manifolds and $F=F_0*F_1$ is a Murasugi sum of $\pi_1$-essential spanning surfaces $F_i\subset M_i$. Write $\min_{i=0,1}\text{ess}(F_i)=n$. Then $F$ is essential, and in fact $\ess(F)\geq n$. Moreover, if neither $F_i$ contains an essential curve that is $\partial$-parallel in $M_i$-cut-along-$F_i$, then the same is true of $F$ in $M$-cut-along-$F$.
\end{T:PlumbEndEss}

In particular, this allows us to extend the main result of \cite{endess} (see \textsection\ref{S:3Mfld} for the relevant definitions):

\begin{T:Endess}
If $P\subset \Sigma$ is a cellular alternating diagram without removably nugatory crossings, then every adequate state surface from $P$ is end-essential in $\Sigma\times I$. 
\end{T:Endess}

Previous versions of this paper included the main result of \cite{natural}, which establishes that a Murasugi sum of two geometrically essential spanning surfaces need not be geometrically essential. When splitting the earlier, combined manuscript into two parts, some overlap in background material was unavoidable. We have attempted to minimize that overlap by streamlining the background \textsection\ref{S:Background} of this paper, while retaining a detailed treatment of the same material in \textsection\textsection2-3 of \cite{natural}.  

With this in mind, the paper is organized as follows. In \textsection\ref{S:Background}, we quickly review several technical aspects of spanning surfaces.  In \textsection\ref{S:Essences}, we define the algebraic essence and geometric essence of a spanning surface. In \textsection\ref{S:Main}, we prove our main result, Theorem \ref{T:PlumbEss}.  In \textsection\ref{S:Twisted}, we determine the geometric essence of a checkerboard surface $F$ for a reduced alternating diagram, introduce twisted plumbing, and use twisted plumbing to prove Theorem \ref{T:CBEss}.  In \textsection\ref{S:3Mfld}, we extend everything to arbitrary 3-manifolds.

\section{Background}\label{S:Background}

We briefly review the following aspects of spanning surfaces, all of which are described in greater detail in \textsection\textsection2-3 of \cite{natural}: algebraic and geometric notions of essential surfaces, Murasugi sums, and techniques involving capping and cutting.

\subsection{Spanning surfaces}


\begin{definition}\label{D:SpanningSurface}
A {\bf spanning surface} $F$ for a link $L$ is a compact surface, orientable or nonorientable, with no closed components which is embedded with $\partial F=L$. An oriented spanning surface is called a {\bf Seifert surface}. We regard two spanning surfaces $F$ and $F'$ as equivalent, denoted $F\approx F'$, if they are {\bf freely isotopic}, meaning that there is an ambient isotopy of $S^3$ that takes $F$ to $F'$.
\end{definition}


\begin{notation}\label{N:Span}
Throughout, given a submanifold $S$ of an some ambient manifold, $\nu S$ will denote a closed regular neighborhood, $\mathring{\nu}S$ will denote the interior of this neighborhood, and $\mathring{S}$ will denote the interior of $S$.  The notations $F$ and $L$ will always denote a spanning surface in $S^3$ and its boundary. In \textsection\textsection\ref{S:Background}-\ref{S:Twisted}, we work in $S^3$.
\end{notation}

\begin{rem}\label{R:Exterior}
It is common to think of a spanning surface $F$ for a link $L\subset S^3$ as being properly embedded in the link exterior $S^3\setminus\mathring{\nu} L$, so that $\partial F$ intersects each meridian on $\partial\nu L$ transversally in one point. We will find it more convenient, however, to take the perspective of Definition \ref{D:SpanningSurface}, except where we specify otherwise.
\end{rem}

One may cut $S^3$ along a spanning surface $F$ to obtain a compact 3-manifold with boundary, which we denote by $S^3\cut F$. Formally, this space is the metric closure of $S^3\setminus F$, where points are equivalence classes of Cauchy sequences that are identified if their interpolation converges. It is homeomorphic to $S^3\setminus\mathring{\nu}F$, but with extra structure from $F$ and $L$ encoded in its boundary. When $F$ is orientable, $S^3\cut F$ is a {\it sutured manifold}, and the extra structure on its boundary is a copy of $L$, which  which cuts $\partial(S^3\cut F)$ into two copies of $F$. When $F$ is nonorientable, however, this copy of $L$ does not separate $\partial (S^3\cut F)$, so $S^3\cut F$ is not quite a sutured manifold.  Nevertheless, this perspective will often prove useful, especially for careful statements of definitions, so we find it worthwhile  to establish the following notation.   

\begin{notation}\label{N:h_F}
Write $ h_F:S^3\cut F\to S^3$ for the quotient map that reglues corresponding pairs of points from $\mathring{F}$ in $\partial (S^3\cut F)$, and denote $\wt{L}={ h_F}^{-1}(L)\subset\partial (S^3\cut F)$ and $\wt{F}= h_F^{-1}(\mathring{F})=\partial (S^3\cut F)\setminus \wt{L}$, so that $h_F$ restricts to a homeomorphism $S^3\cut F\setminus\wt{F}\to S^3\setminus\mathring{F}$ and to a 2:1 covering map $\wt{F}\to\mathring{F}$.
\end{notation}

\begin{figure}[h]
\begin{center}
{\includegraphics[height=1in]{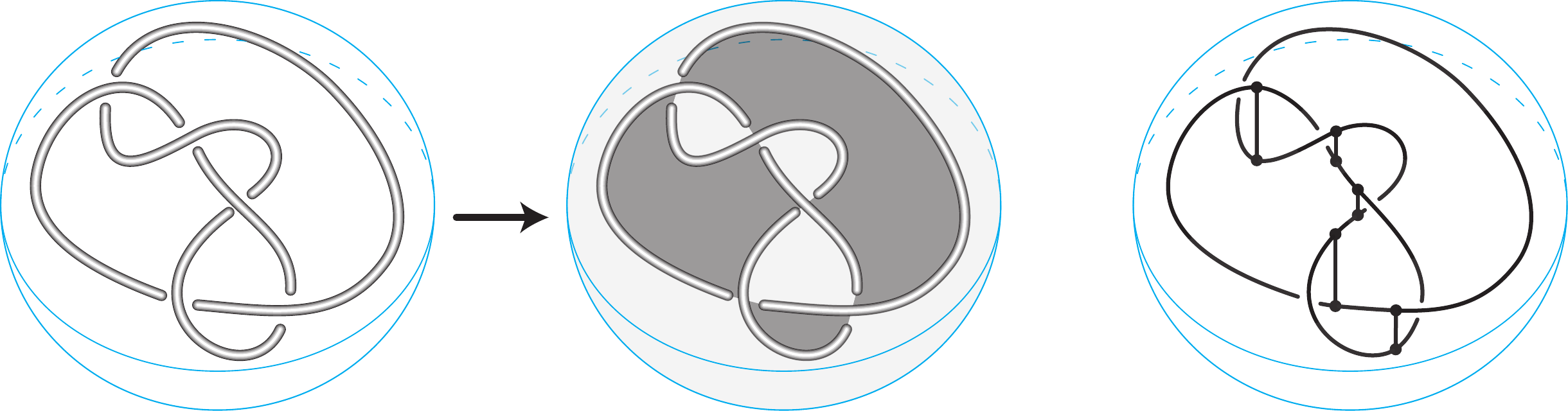}}
\caption{Constructing checkerboard surfaces}
\label{Fi:Chessboard}
\end{center}
\end{figure}

\begin{definition}
Given a diagram $P$ of $L$, one can construct two spanning surfaces $B$ and $W$ by coloring the regions of $S^2\setminus P$ black and white in checkerboard fashion. The interiors of these {\bf checkerboard surfaces} $B$ and $W$ intersect in {\it vertical arcs} which project to the crossings of $P$. 
\end{definition}

Figure \ref{Fi:Chessboard} shows this construction and the spatial graph $B\cap W$ comprised of $L$ and the vertical arcs at the crossings. The construction generalizes as follows:

\begin{definition}
Given a diagram $P$ of a link $L$, smoothing each crossing  in one of two ways, $\raisebox{-.02in}{\includegraphics[width=.125in]{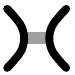}}
\overset{\color{Gray}{_{{A}}}\color{black}}{\longleftarrow}\raisebox{-.02in}{\includegraphics[width=.125in]{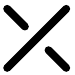}}
\overset{_{{B}}}{\longrightarrow}\raisebox{-.02in}{\includegraphics[width=.125in]{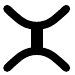}}$
yields a {\it state} $x$ of $P$, comprised of {\it state circles} and \color{Gray}$A$\color{black}- and $B$-labeled edges. One constructs an associated {\it state surface}  $F_x$ by capping off the state circles with mutually disjoint {\it state disks} (that are also disjoint from crossings and transverse to the projection sphere) and attaching a half-twisted band at each crossing \cite{ozawa11,ak}. %
The isotopy class of $F_x$ may depend on the {\it layering} of the disks relative to the projection sphere; to avoid such ambiguity, we assume, unless stated otherwise, that all state circles are capped with disks {\it lying entirely on the same side} of the projection sphere $S^2$. 
\end{definition}

\begin{prop}[Proposition 1.3.1 of \cite{TkThesis}, Proposition 2.6 of \cite{natural}]\label{P:StateToCB}
Any state surface $F_x$ from any diagram $P$ is (freely) isotopic to a checkerboard surface of some diagram $P'$.
\end{prop}

See \cite{natural} or \cite{TkThesis} for a simple proof that constructs the needed isotopy.

\begin{definition}\label{D:Homogeneously_Adequate}
Given a state $x$ of a link diagram
, the (abstract) {\bf state graph} $\Gamma_x$ is obtained by collapsing each state circle to a point, while keeping the ${\color{myGray} A}$ and $B$ labels on the edges.
The state $x$ is {\bf adequate} if $\Gamma_x$ has no loops (i.e. the endpoints of each $A$- and $B$-labeled edge lie on distinct state circles), and $x$ is {\bf homogeneous} if all edges in each cut component\footnote{If $\Gamma_x$ has no cut vertices
, then $\Gamma_x$ has a single cut component; otherwise, cut $\Gamma_x$ at a cut vertex. Cut each resulting component at a cut vertex, if one exists.  Continue until no component has a cut vertex.  The resulting components are the {\it cut components} of $\Gamma_x$.} of $\Gamma_x$ have the same type, ${\color{myGray} A}$ or $B$.  If both conditions hold, $x$ and $F_x$ are called {\bf homogeneously adequate}; in this case, $F_x$ is $\pi_1$-essential \cite{ozawa11}.
\end{definition}

\subsection{Geometrically and algebraically essential surfaces}

There are two common notions of essentiality  for properly embedded surfaces in a 3-manifold. For further discussion and motivation, see \textsection1 of \cite{ht}, \textsection2 of \cite{ak}, or \textsection2.2 and the Appendix of \cite{natural}.

\begin{figure}[h]
\begin{center}
\includegraphics[height=.6in]{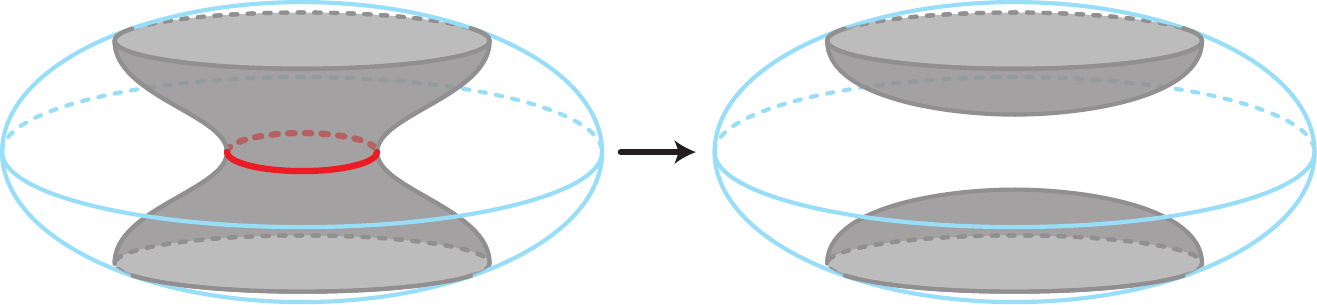}\\
\includegraphics[height=.6in]{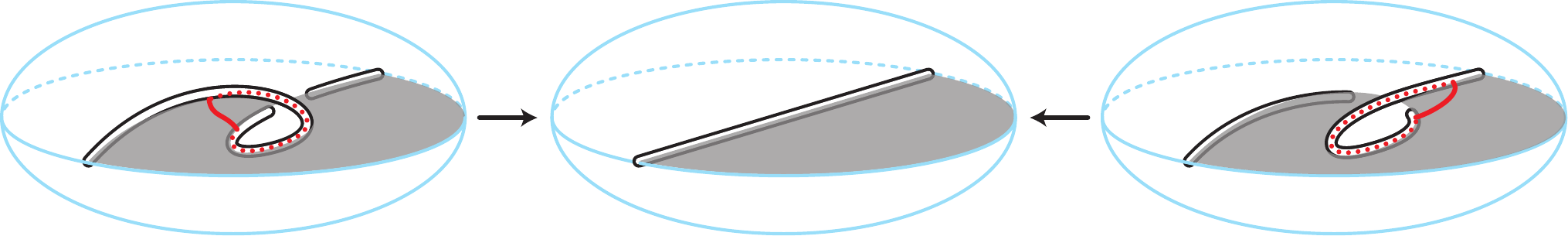}
\caption{Geometric compression (top) and $\partial$-compression (bottom) of a spanning surface}
\label{Fi:Compressions}
\end{center}
\end{figure}

\begin{definition}\label{D:GeomEss}
A surface $F$ spanning a link $L$ is {\bf geometrically essential} if $F$ cannot be compressed or ${\partial}$-compressed to a spanning surface (see Figure \ref{Fi:Compressions}).
That is, $F$ is geometrically essential if {\it both}:
\begin{enumerate}
\item For every embedded disk $D\subset S^3\setminus L$ with $D\cap F=\partial D$, the circle $\partial D$ bounds a disk in $F$, and\footnote{We use ``circle" as shorthand for ``simple closed curve". A circle in a surface is {\it essential} if it does not bound a disk in that surface.}
\item For every embedded disk $D\subset S^3\setminus L$ with $\partial D=\alpha\cup\beta$ for arcs $\alpha\subset F$ and $\beta\subset L$, the arc $\alpha$ is $\partial$-parallel in $F$.
\end{enumerate}
If $F$ satisfies (1), it is called {\bf geometrically incompressible}, whether or not it satisfies (2). Otherwise, there is an essential circle in $F$ that bounds a disk $D\subset S^3\setminus L$ with $D\cap F=\partial D$, called a (geometric) {\bf compressing disk} for $F$.
\end{definition}

\begin{rem}
If $F$ is geometrically incompressible but geometrically inessential, then  $F$ is isotopic to $F'\natural\MobPos$ or $F'\natural\MobNeg$ for some spanning surface $F'$.
\end{rem}

\begin{definition}\label{D:AlgEss}
A surface $F$ spanning a link $L$ is {\bf $\boldsymbol{\pi_1}$-essential} if :
\begin{enumerate}
\item Inclusion ${\mathring{F}\hookrightarrow S^3\setminus L}$ induces an injection of fundamental groups, and 
\item $F$ is not a M\"obius band spanning the unknot, $F\not\approx\MobPos,\MobNeg$.
\end{enumerate}
If $F$ satisfies (1), it is called {\bf $\pi_1$-injective}, whether or not it satisfies (2).
\end{definition}

\begin{rem}
If $F$ is ${\pi_1}$-essential, {then} $F$ is geometrically essential.  The converse holds when $F$ is orientable by the loop theorem, but not in general \cite[Prop. 4.2]{natural}.
\end{rem}

\begin{rem}
The loop theorem implies that, given a $\pi_1$-non-injective spanning surface $F$, there is a smooth map $f:D^2\to S^3\setminus L$ which restricts to an embedding $f:\mathring{D}^2\hookrightarrow S^3\setminus F$ of the disk's interior and to an immersion $f:S^1\to F$ of the disk's boundary $S^1=\partial D^2$.  Given such a map $f$, we call its image $D=f(D^2)$ an {\bf algebraic compressing disk}. It is just like a typical compressing disk, except that $\partial D$ is immersed in $\mathring{F}$ and may (or may not) self-intersect.
\end{rem}

\subsection{Murasugi sum}\label{S:Plumb}

\begin{figure}[h]
\begin{center}
\labellist
\pinlabel{$=$} at 165 30
\pinlabel{$*$} at 325 30
\endlabellist
\includegraphics[width=.8\textwidth]{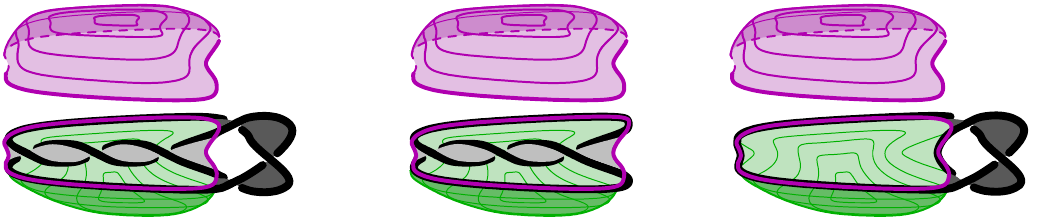}
\caption{Decomposing under Murasugi sum, or deplumbing}
\label{Fi:Plumb}
\end{center}
\end{figure}

\begin{definition}\label{D:Plumb}
Let $V\subset S^3$ be an embedded disk with $V\cap F=\partial V$ such that
\begin{enumerate}
\item $\partial V$ bounds a disk $U\subset F$.
\item Denoting $S^3\cut(U\cup V)=B_0\sqcup B_1$, neither $F_i=F\cap B_i$ is a disk.
\end{enumerate}
Then $U$ is a {{\bf plumbing disk}} for $F$, and $V$ is an associated {{\bf plumbing cap}}. If $V$ satisfies (1) but not (2), we call $V$ a {\it fake plumbing cap}.

Say that $F$ is the {\bf Murasugi sum} of $F_0$ and $F_1$ along $U$, and write $F_0*F_1=F$.  This operation 
is also called {\bf generalized plumbing}, although it is often convenient for linguistic reasons to use various forms of the phrasing ``plumbing'' more liberally, provided this does not lead to confusion. For example, we may call the associated decomposition a {\bf deplumbing} (see Figure \ref{Fi:Plumb}).\end{definition}

\begin{rem}\label{R:CBPlumb}
Let $x$ be a state, $F_x$ its state surface, and $\Gamma_x$ its state graph.  Each cut component $\Gamma_i$ has an associated subsurface $F_i\subset F_x$.  The surface $F_x$ can be constructed by plumbing the subsurfaces $F_i$ in such a way that, if $\Gamma_i\cap\Gamma_j=\{u\}$, then $F_i$ and $F_j$ are plumbed along the state disk $U$ that corresponds to the vertex $u\in \Gamma_x$.  The key point here is that, for every non-innermost state disk $U\subset F_x$, there is another disk $V$ on the opposite side of the projection sphere such that $F_x$ deplumbs in a nontrivial way along $U\cup V$.  (If one allows the state disks to intersect $S^2$, the same argument still works, taking $V$ to lie entirely on one side or the either of $S^2$.) See Figure \ref{Fi:Plumb}.
\end{rem}

A great deal is known about Murasugi sums, especially oriented ones. See \cite{natural} for a survey.  Here, we highlight an important result of Gabai, along with an extension of this result and a corollary thereof by Ozawa:

\begin{theorem}\label{T:gabai}[\cite{gab1,gab2}]
If $F_0*F_1=F$ is a Murasugi sum of Seifert surfaces with each $\partial F_i=L_i$ and $\partial F=L$, then:
\begin{enumerate}[label=(\arabic*)]
\item $F$ is {essential} if $F_0$ and $F_1$ are essential.
\item $F$ has {minimal} genus if and only if $F_0$ and $F_1$ both have minimal genus.
\item $L$ is a {fibered link with fiber} $F$ {if and only if} each $L_i$ is fibered with fiber $F_i$.
\item $S^3\setminus \overset{_\circ}{\nu}L$ has a {nice codimension 1 foliation}  {if and only if} both $S^3\setminus \overset{_\circ}{\nu}L_i$ do.
\end{enumerate}
\end{theorem}

\begin{theorem}
[Lemma 3.4 of \cite{ozawa11}]\label{T:Ozawa}
If $F=F_0*F_1$ is a Murasugi sum of $\pi_1$-essential spanning surfaces $F_i$, then $F$ is $\pi_1$-essential.
\end{theorem}

\begin{theorem}[Theorem 2.8 of \cite{ozawa11}]\label{T:ozawafkp}
If $x$ is a homogeneously adequate state, then the state surface $F_x$ (with any layering of the state disks) is $\pi_1$-essential.
\end{theorem}

\subsection{Caps and height}\label{S:Caps}

Next, we introduce a few related technical ideas that we will need for our definitions and proofs.  The underlying ideas here are traditional.  An interested reader can find a more extensive application of these techniques in the paper \cite{natural}, which was previously connected to this one. A formal treatment of these techniques appears in still greater detail in the author's doctoral thesis \cite{TkThesis}.

\begin{definition}\label{D:cap}
A {\bf cap} for $F$ is the image $V= h_F(\wt{V})$ of a compressing disk for $\partial (S^3\cut F)$. We always assume that $\partial\wt{V}\pitchfork\wt{L}$ in $\partial (S^3\cut F)$.
\end{definition}

\begin{figure}[h]
\begin{center}
\includegraphics[height=.5in]{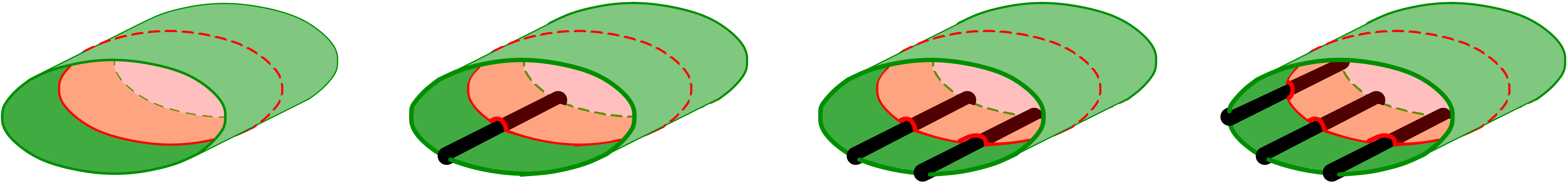}
\caption{Compressing disks for $\partial (S^3\cut F)$, whose images under $h_F$ are called caps for $F$}
\label{Fi:AlgCaps}
\end{center}
\end{figure}

See Figure \ref{Fi:AlgCaps}. Note that if $V$ is a cap for $F$, then we allow $\partial V$ to intersect itself or $L$. Note also that if $\partial V\cap L=\varnothing$, then $\partial V$ cannot be contractible in $F$, or else $\partial\wt{V}$ would be contractible in $\partial (S^3\cut F)$.  If $\partial V$ intersects $L$, however, then $\partial V$ may well be contractible in $F$.  For example, this is the case if $V$ is a plumbing cap.

\begin{notation}\label{N:wt(D)}
If $D$ is a cap for $F$, then $\wt{D}$ denotes the (unique) properly embedded disk in $S^3\cut F$ satisfying $ h_F(\wt{D})=D$.
\end{notation}

\begin{figure}[h]
\begin{center}
\labellist \tiny \hair 4pt
\pinlabel {$0$} at 50 100
\pinlabel {$\White{3}$} at 200 95
\pinlabel {$2$} at 252 40
\pinlabel {${2}$} at 330 78
\pinlabel {$\White{1}$} at 162 30
\pinlabel {$0$} at 127 12
\pinlabel {$0$} at 200 12
\pinlabel {$\White{0}$} at 290 14
\pinlabel {$\White{1}$} at 330 103
\pinlabel {$0$} at 323 126
%
\endlabellist
\includegraphics[height=1in]{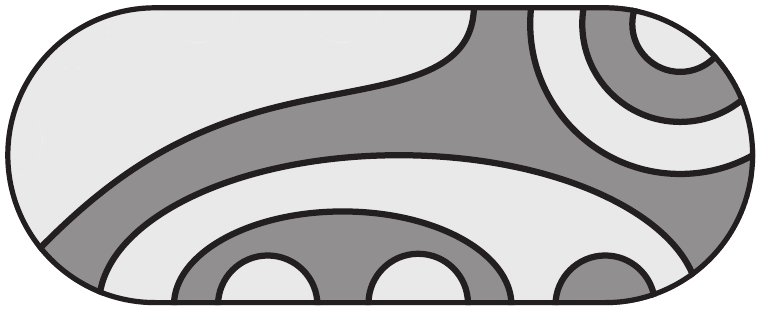}\hfill
\caption{A disk $D$ cut by arcs and labeled by height}
\label{Fi:TreeLabel}
\end{center}
\end{figure}

\begin{definition}\label{D:height}
Given a disjoint union $A=\bigsqcup_{i\in I}\alpha_i$ of $\ell$ properly embedded arcs in a disk $D$, let $T$ be the tree with one vertex for each subdisk comprising $D\cut A$ in which two vertices are adjacent whenever the corresponding subdisks abut.  

Define the {\bf height} of each subdisk of $D\cut A$ recursively as follows. Let $T_0=T$. Outermost disks of $D\cut A$, corresponding to leaves in $T_0$, have height 0.   For $i\geq 1$, let $T_{i}$ be the tree obtained from $T_{i-1}$ by deleting each leaf and its incident edge.  Disks of $D\cut A$ that correspond to leaves in $T_i$ have height $i$.   Figure \ref{Fi:TreeLabel} shows an example.
\end{definition}


\section{Algebraic and geometric essence}\label{S:Essences}
 
\subsection{Types of caps}
Recall from Notations \ref{N:h_F} and \ref{N:wt(D)} and Definitions \ref{D:Plumb} and \ref{D:cap} that a {\it cap} $D$ is the image under the natural map $h_F:S^3\cut F\to S^3$ of a compressing disk for $\partial (S^3\cut F)$, that $\wt{D}$ denotes the (unique) compressing disk for $\partial(S^3\cut F)$ satisfying $h_F(\wt{D})=D$, and that $D$ is a {\it plumbing cap} if $\partial D$ bounds a disk in $F$. 
More generally, define the following types of caps:

\begin{definition}\label{D:Caps}
Let $D$ be a cap for $F$.  Call $D$ {\bf $\boldsymbol{\partial}$-contractible} if $\partial D$ is contractible in $F$, {\bf $\boldsymbol{\partial}$-essential} if it does not, and {\bf geometric} if $\partial D$ does not self-intersect.
\end{definition}

\begin{figure}[h]
\begin{center}
\includegraphics[width=\textwidth]{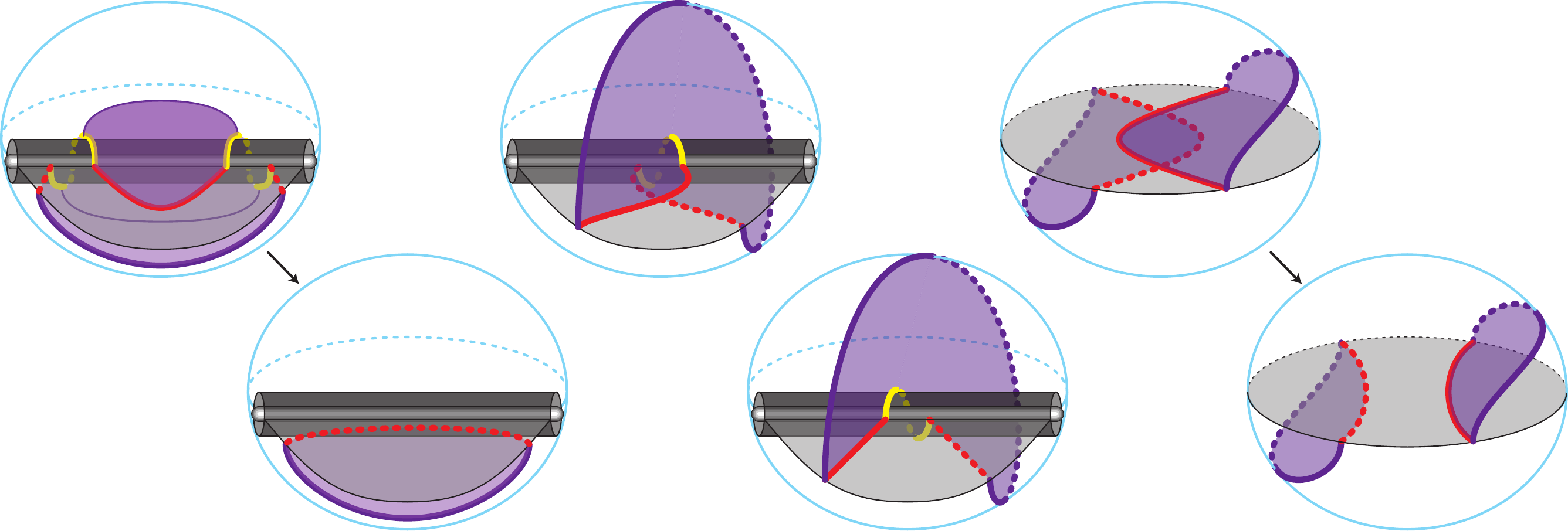}
\caption{A cap is {\it acceptable} if it cannot be simplified by these moves.}
\label{Fi:Acceptable}
\end{center}
\end{figure}

\begin{definition}\label{D:Acceptable}
Let $D$ be a cap for $F$, 
denote the set of self-intersection points of $\partial D$ by $x$, and take $\nu L$ to be generic.\footnote{This implies that no arc of $D\cap\partial \nu L$ is parallel in $\partial \nu L$ to $F\cap\partial \nu L$, because we assume in Definition \ref{D:cap} that $\partial\wt{V}\pitchfork\wt{L}$ in $\partial (S^3\cut F)$.}  We say that $D$ is {\bf acceptable} if $D$ admits none of the simplifying moves shown in Figure \ref{Fi:Acceptable}: 
\begin{itemize}
\item no arc of $\partial D\setminus\mathring{\nu}L$ cuts off a disk from $F\setminus\mathring{\nu} L$ (see Figure \ref{Fi:Acceptable}, left); and
\item no two arcs of $\partial D\cut (x\cup \nu L)$ together cut off a disk from $F\setminus\mathring{\nu} L$ (center and right).
\end{itemize}
\end{definition}


\begin{notation}\label{N:Caps}
Denote the set of all caps for $F$ by $\Caps(F)$. Also denote the set of geometric caps for $F$ by $\Caps_g(F)$ and the set of $\partial$-essential caps by $\Caps_e(F)$.
\end{notation}

\subsection{Fake caps}

Recall that a {\it fake plumbing cap} is a $\partial$-contractible geometric cap $D$ for $F$ such that $F$ intersects one of the closed balls comprising $S^3\cut (D\cup U)$ in a disk. In this case, the disk $\partial \wt{D}$ (obtained by lifting $D$ via $h_F$) bounds a disk in $\partial (S^3\cut F)$. Generalize this terminology as follows.  

\begin{definition}
If $\wt{D}$ is a properly embedded disk in $S^3\cut F$ such that $\partial\wt{D}$ bounds a disk in $\partial (S^3\cut F)$, then $D= h_F(\wt{D})$ is a {\bf fake cap} for $F$.
\end{definition}

\begin{obs}\label{O:FakeDegree}
If $D$ is a fake cap for $F$, then $D$ is $\partial$-contractible and $|\partial D\cap L|$ is even.
\end{obs}

Note that if $D$ is a fake cap for $F$ with $\partial D\cap L= \varnothing$ and $L$ is non-split, then $D$ is parallel through $S^3\cut F$ into $\mathring{F}$. On the other hand, we find it instructive to note the following.

\begin{prop}\label{P:FakeAcc}
Suppose $D$ is a fake cap for $F$. If $\partial D\cap L\neq \varnothing$, then $D$ admits a simplifying move from Figure \ref{Fi:Acceptable}.
\end{prop}

That is, ``acceptable fake caps'' do not exist. Since we will not use this fact, we leave the proof to the reader. 

We will use the contrapositive of the following proposition when we prove Theorem \ref{T:TwistedEss}.

\begin{prop}\label{P:FakeString}
Suppose $D$ is a fake cap for $F$ with $|\partial D\cap L|=2r$, so that $\wt{L}$ cuts $\partial \wt{D}$ into $2r$ arcs. Then any constructive string of $r$ of these arcs of $\partial \wt{D}\cut \wt{L}$ contains at least one arc whose image under $h_F$ is $\partial$-parallel in $F$. 
\end{prop}

\begin{proof}
This is vacuous when $r=0$. Otherwise, such a string of $r$ arcs is incident to $r+1$ endpoints of arcs of $\wt{L}\cap\wt{D}$.  Since there are just $r$ such arcs, at least one must have both its endpoints on the string in question. Therefore, some {\it outermost} (in $\wt{D}$) arc of $\wt{L}\cap\wt{D}$ has both endpoints on this string. It is parallel through a disk $\wt{D}_0$ to one of the $r$ arcs in our string.  The image of that arc under $h_F$ will be parallel through $h_F(D_0)$ into $\partial F$.
\end{proof}

\subsection{Essence}\label{S:Essence}
 
 Given a cap $D$ for $F$, we call $|\partial D\cap L|$ the {\bf complexity} of $D$. It is sensible to consider caps of minimal complexity:

 \begin{definition}\label{D:Essence}
The (algebraic) {\bf essence} of $F$ is
\[\text{ess}(F)=\min_{D\in\Caps_e(F)}|\partial D\cap L|.\]
Likewise, the {\bf geometric essence} of $F$ is 
\[\text{ess}_g(F)=\min_{D\in\Caps_g(F)}|\partial D\cap L|,\]
and the {\bf $\boldsymbol\partial$-contractible essence of $F$} is
\[\text{ess}_c(F)=\min_{D\in\Caps(F)\setminus\Caps_e(F)}|\partial D\cap L|.\]
\end{definition}
We will show in Theorem \ref{T:PlumbEss} that $\ess(F)$ is well-behaved under plumbing.  
First, several remarks:

\begin{rem}
Every spanning surface $F$ satisfies $\text{ess}(F)\leq \text{ess}_g(F)$ because $\Caps(F)\subset \Caps_g(F)$.
\end{rem}


\begin{rem} 
A spanning surface $F$ is $\pi_1$-essential if and only if $\text{ess}(F)\geq 2$, and $\text{ess}(F)=1$ if and only if $F$ is a M\"obius band spanning the unknot, $F=\MobPos,\MobNeg$.
 \end{rem}

\begin{rem}
A spanning surface $F$ is geometrically incompressible if and only if  $\text{ess}_g(F)\geq 1$ and is geometrically essential if and only if $\text{ess}_g(F)\geq 2$.
\end{rem}

Theorem 4.1 of \cite{natural} implies that plumbing does not respect geometric essence.  Yet, geometric essence has some advantages, relative to (algebraic) essence. One is that any acceptable geometric cap $D$ describes a possible surgery move on $F$, generalizing the way that geometric compressing disks and $\partial$-compressing disks do.

\begin{rem}
Viewing $F$ as a properly embedded surface in the link exterior $S^3\setminus\mathring{\nu}L$ (just for this paragraph), one can cut $F$ along the $n$ arcs of $\partial D\cap F$ and glue in two parallel copies of $D$.  
The resulting surface $F'$ satisfies $\beta_1(F')=\beta_1(F)+n-2$, and for each component $L_i$ of $L$, the component-wise boundary slopes $s(F,L_i)=\text{lk}(L_i,F\cap\partial\nu L_i)$ and $s(F',L_i)=\text{lk}(L_i,F'\cap\partial\nu L_i)$ satisfy $|s(F,L_i)-s(F',L_i)|\leq 2|\partial D\cap L_i|$. 
A much more detailed discussion of these surgery moves will appear in \cite{cklsv3}, using results from \cite{cklsv2}. 
\end{rem}

Another advantage of geometric, over algebraic, essence, is its relationship to framings, Gordon-Litherland matrices, and (in the case of state surfaces) state graphs.  For background on these topics, see \cite{greene} or \cite{flyping}, e.g. We discuss these relationships further in \textsection\ref{S:Twisted}.

\section{Plumbing respects essence}\label{S:Main}

Now we 
prove the following generalization of Ozawa's essential plumbing theorem \cite{ozawa11}. 
%
We note two reasons why the proof is somewhat difficult.  The first reason is that it concerns a plumbing $F=F_0*F_1$ along a sphere $U\cup V$, where $|\partial V\cap L|$ can be arbitrarily large, independent of the essence of either surface.  This is why, when considering a cap $D$ for $F$ of minimal complexity and the subdisks into which $D$ cuts $V$, it will not suffice to consider outermost subdisks of $V\cut D$.  Instead, a somewhat more subtle counting argument is needed.  The same counting method underpins Ozawa's proof.  The second reason is that we will encounter fake caps in a way that might seem obviously contradictory, but distilling contradictions from them will take some extra work.  Indeed, other than careful setup and this extra work regarding fake caps, the proof takes just a few sentences.

\begin{theorem}\label{T:PlumbEss}
Suppose $F=F_0*F_1$ is a plumbing of $\pi_1$-essential spanning surfaces, and write $\min_{i=0,1}\text{ess}(F_i)=n$. Then $\ess(F)\geq n$.
\end{theorem}

\begin{proof}
Let $U$ and be a plumbing disk along which $F$ decomposes into $F_0$ and $F_1$, let $V$ be an associated plumbing cap, and write $S^3\cut(U\cup V)=B_0\sqcup B_1$ such that each $F_i=F\cap B_i$, and let $S^3_i$ denote a copy of $S^3$ containing just $F_i$. 
We use $U_i$ to denote the copy of $U$ in $\partial B_i\cut F_i$. See Figure \ref{Fi:1D_Murasugi_Sum}.
    
    \begin{figure}[h]
        \centering
        \labellist \small \hair 4pt
        \pinlabel {$=$} at 315 470
        \pinlabel {$*$} at 665 470
        \pinlabel {$B_0$} at 140 530
        \pinlabel {$B_1$} at 130 320
        \pinlabel {$F_0$} at 460 420
        \pinlabel {$F_1$} at 870 310
        \pinlabel {{\color{myGreen} $U$}} at 140 390
        \pinlabel {{\color{myBlue} $V$}} at 250 480
        \pinlabel {{\color{myGreen} $U_0$}} at 140 145
        \pinlabel {{\color{myGreen} $U_1$}} at 140 90
        \endlabellist
        \includegraphics[width=.6\textwidth]{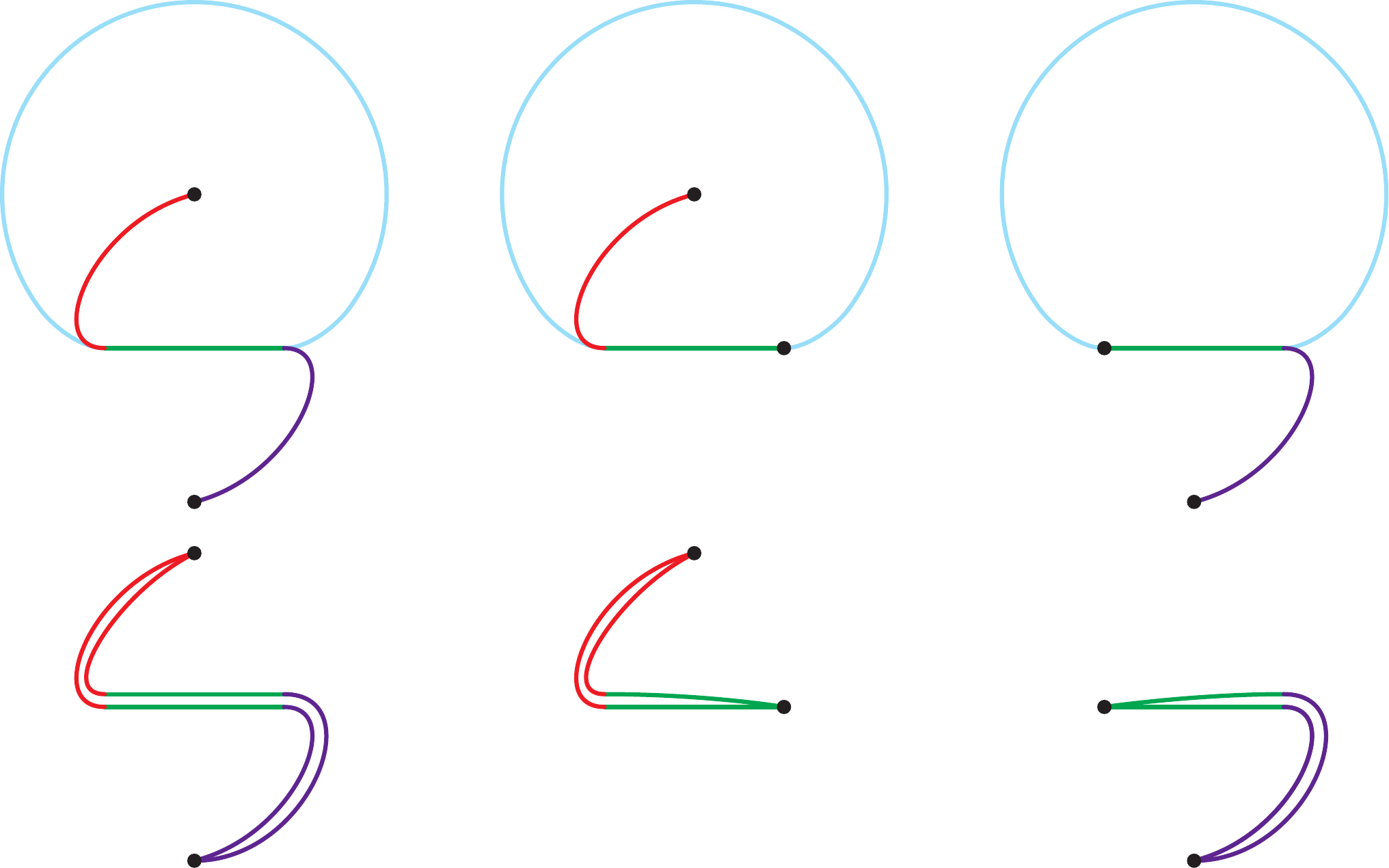}
        \caption{A schematic picture of a Murasugi sum in $S^3$, labeled for the proof of Theorem  \ref{T:PlumbEss}}
        \label{Fi:1D_Murasugi_Sum}
    \end{figure}


Choose $D\in \Caps_e(F)$ so as lexicographically to minimize $k=|\partial D\cap L|$, $\ell=|D\cap V|$, and $m=|\partial D\cap\partial U|$, provided that $D\pitchfork V$, and that $\partial U$ contains no points where $\partial D$ intersects itself or $L$. 
 If $\ell=0$, then \textsc{wlog} $D$ is a cap for $F_1$ and $\partial D$ is not contractible in $F_1$, so $\ess(F)= |\partial D\cap L|\geq \ess(F_1)\geq n$, and we are done.  Assume instead that $\ell>0$. Denote $\partial D=\gamma$.

Denote the $\ell$ arcs of $D\cap V$ by $\alpha_1,\dots,\alpha_\ell$.  Choose (say, mutually disjoint) arcs $\beta_1,\dots,\beta_\ell\subset U$ such that $\partial \beta_i=\partial\alpha_i$.

\begin{figure}[h]
\begin{center}
\includegraphics[width=.8\textwidth]{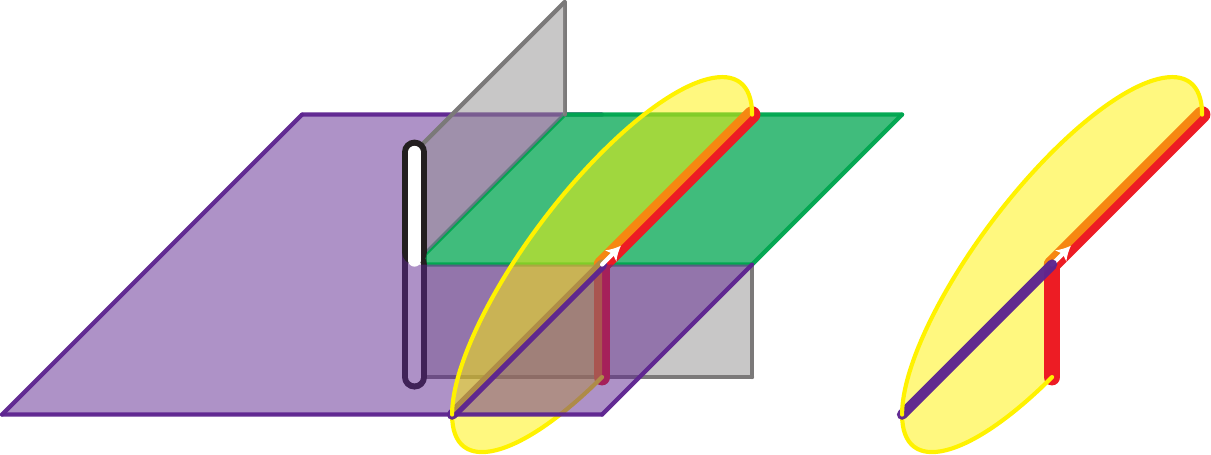}
\caption{
Given a cap $D$ (yellow), mark ${\color{myRed} \partial D}$ (red) near each endpoint of ${\color{myPurple} D\cap V}$ (purple; all of $V$ is purple) with an arrow (white) that points in the direction where ${\color{myRed} \partial D}$ runs along ${\color{myGreen} U}$ (green). 
}
\label{Fi:Marker}
\end{center}
\end{figure}

 The arcs $\alpha_1,\dots,\alpha_\ell$ cut $D$ 
into $\ell+1$ disks, $D_0,\dots,D_\ell$.  Decorate each of the $2\ell$ endpoints of $V\cap D$ on $\gamma$ with an arrow (``marker") that points in the direction where $\gamma$ runs along $U$, as in Figure \ref{Fi:Marker}). 
Add one additional marker in each of the $k$ components of $\gamma\cap L$.  There are now $2\ell+k$ or markers, each of which lies on the boundary of exactly one of the disks $D_0,\dots,D_\ell$. See Figure \ref{Fi:OM}.

    \begin{figure}[h]
    \centering
    \includegraphics[height=1.25in]{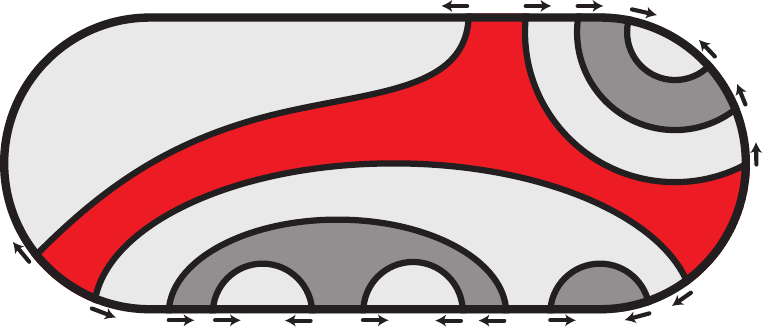}
    \caption{An example of the disk $D$, its arcs of intersection with the disk $V$, and markers at the endpoints of these arcs. Here, just one component of $D\cut V$, colored red, has fewer than two markers.}
    \label{Fi:OM}
\end{figure}

Consider any one of these disks $D_i$, and assume without loss of generality that $D_i\subset B_0$.  See Figure \ref{Fi:D_0}. Each arc $\alpha_t$ of $D\cap V$ on $\partial D_i$ is parallel through a properly embedded disk (``bigon'') in $B_1\subset S^3_0$ (i.e. in $B_1$ ignoring $F\cap\mathring{B}_1$) to the arc $\beta_t\subset U$ (which it may be helpful to think of as lying on $U_1$). The union of $D_i$ and all of these bigons is a cap or fake cap $E_i$ for $F_0$. Moreover, $|\partial E_i\cap \partial F_0|$ equals the number of markers on the boundary of $D_i$.

    \begin{figure}[h]
        \centering
        \labellist
        \pinlabel{$B_0$} at -50 700
        \pinlabel{$B_1$} at -50 0
        \endlabellist
        \includegraphics[width=\textwidth]{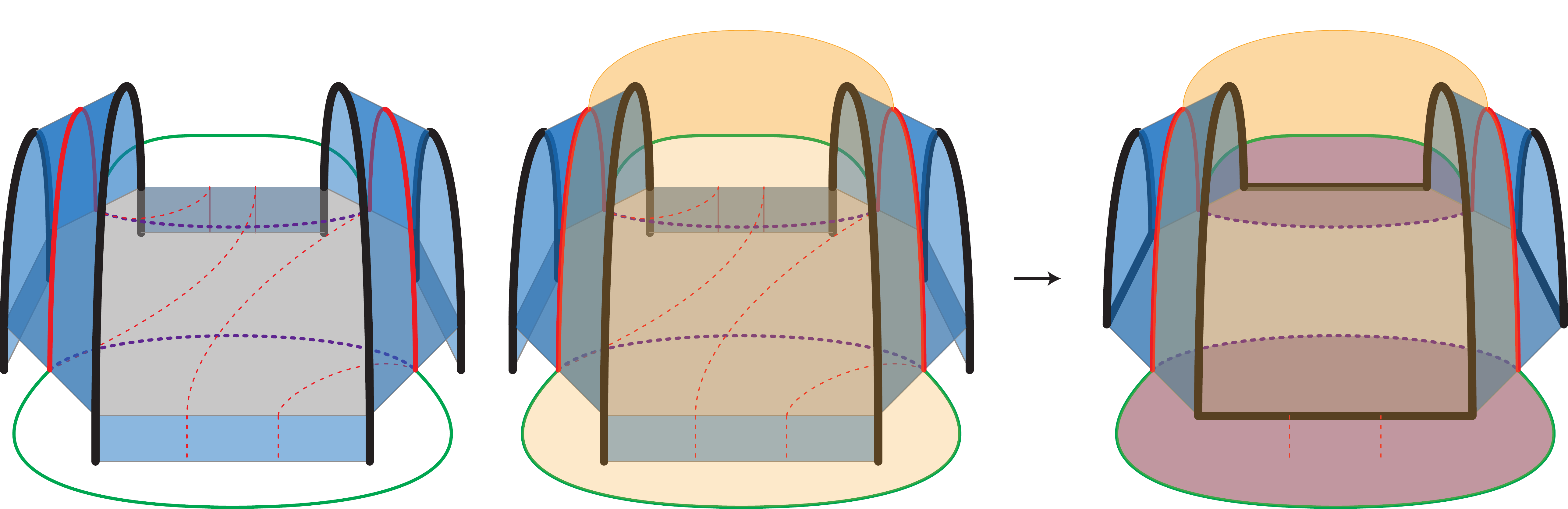}
        \caption{Extending a disk $D_i$ (orange) of $D\cap B_0$ to obtain a compressing disk $E_i$ for $F_0$: each arc of $\partial D_i\cap V$ (green) is parallel through a disk (purple) in $B_1$ to an arc (purple) in $U$; $\partial D$ is red, and $U\cup V$ is a horizontal plane with $U$ gray 
and $V$ not pictured.}
        \label{Fi:D_0}
    \end{figure}
    
Note that the cap or fake cap $E_i$ is not necessarily geometric, even if $D$ is a geometric cap for $F$, because the arcs of $\partial E_i\cap U_1$ might be forced to intersect arcs of $\partial D_i\cap U_0$. This is why the this theorem, like Ozawa's Theorem \ref{T:Ozawa}, holds in the algebraic setting but not the geometric.

We claim that the boundary of each of the disks $D_0,\dots,D_\ell$ contains at least two markers.  To see why, suppose that there are fewer than two markers on $\partial D_i$. It follows that $E_i$ must be a fake cap for $F_0$, since $F_0$ is $\pi_1$-essential. Observation \ref{O:FakeDegree} thus implies that $|\partial E_i\cap L_0|$ is even, hence equals 0 rather than 1.

\begin{figure}[h]
    \centering
        \labellist\small
        \pinlabel{$B_0$} at 330 150
        \pinlabel{$B_1$} at 330 25
        \endlabellist    
        \includegraphics[width=\textwidth]{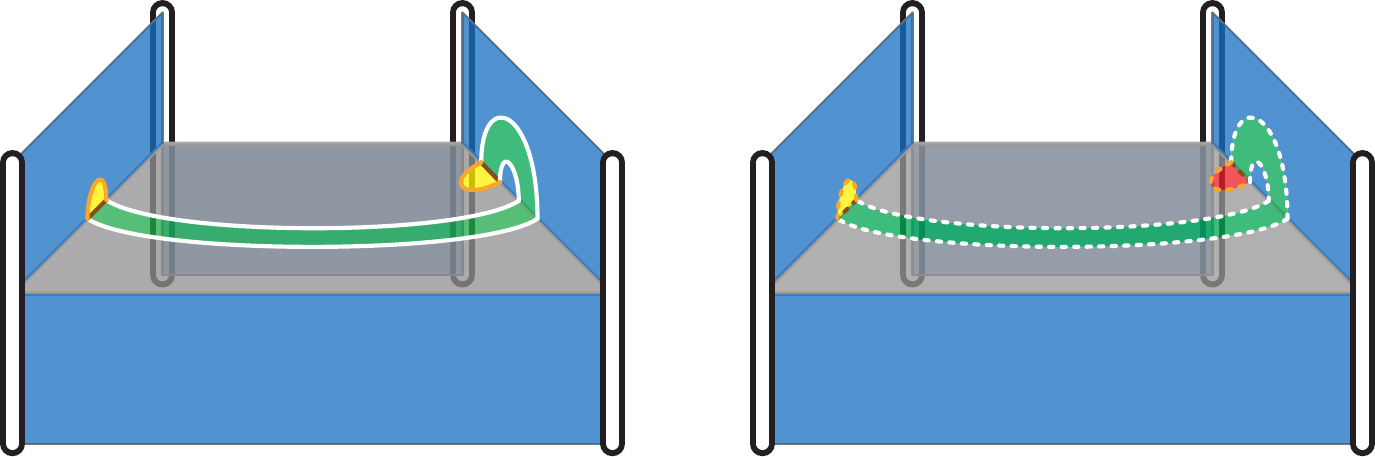}
    \caption{Two simple possibilities for the disk $f(D^2)\subset F_0$, shown green, yellow, and red; an outermost disk $X$ appears yellow or red; the arc $a$ is orange, while $b$ is brown.}
    \label{Fi:Trivial_Disk_F}
\end{figure}

Therefore, there is a continuous map $f:D^2\to F_0$ with $f(S^1)=\partial E_i$. 
Choose $f$ so as to minimize $|f^{-1}(\partial U)|$.  Then $f^{-1}(\partial U)$ is comprised only of arcs, no circles. Figure \ref{Fi:Trivial_Disk_F} shows two simple examples of how this might appear.  Consider an outermost disk of $D^2\cut f^{-1}(\partial U)$, and let $X$ denote the image of this disk under $f$.  This disk $X$ lies either in $U$ or in $F_0\setminus\mathring{U}$, on one side of the surface (locally near $U$) or the other; each of the four possibilities appears as the yellow or red shaded part of one of the two examples in Figure \ref{Fi:Trivial_Disk_F}.

The boundary of $X$ is comprised of two arcs $a$ and $b$, where $a$ is a component of $\partial E_i\cut \partial U$ and $b\subset \partial U\cut \partial E_i$.  We cannot have $a\subset \partial D$ due to the minimality of $|D\cap V|$ and $|\partial D\cap\partial U|$, so $a$ 
    and $X$ must lie in $U_1$ as shown shaded red in Figure \ref{Fi:Trivial_Disk_F}, with $b$ being parallel through one of the bigons of $E_i\setminus\mathring{D}_i$ to an arc component $\alpha_t$ of $D\cap V$. Hence, the circle $\alpha_t\cup b$ bounds a disk $X'$ in $V$.

\begin{figure}[h]
    \centering
    \labellist\small
        \pinlabel{$B_1$} at 330 150
        \pinlabel{$B_0$} at 330 25
        \endlabellist 
    \includegraphics[width=\textwidth]{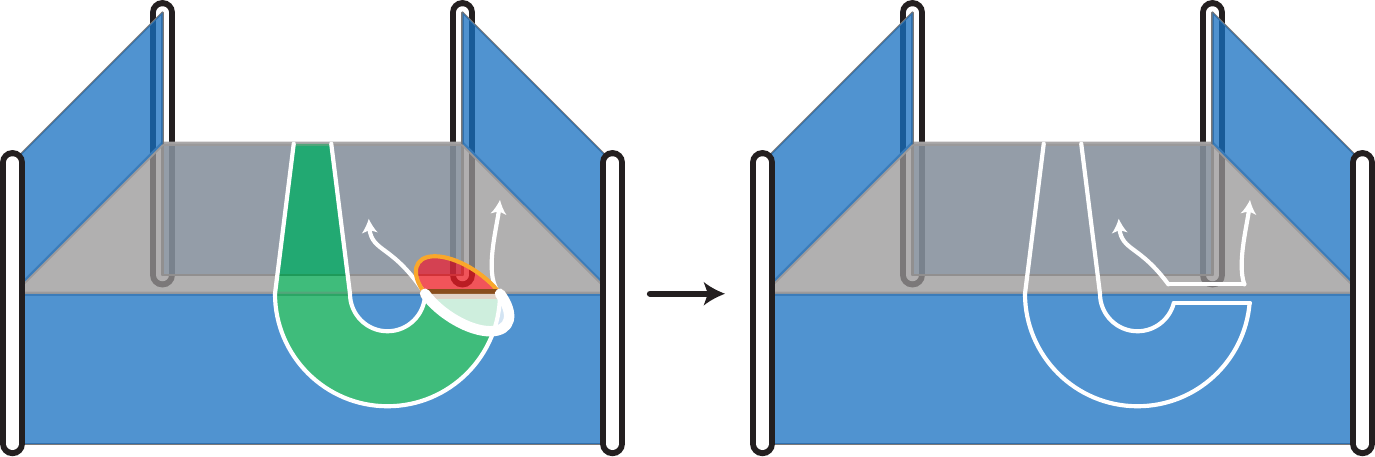}
    \caption{Surgering $D$ along $X'$---the color scheme is the same like Figure \ref{Fi:Trivial_Disk_F}, but the arc $\alpha_t$ (thick), the disk $X'$ (co-bounded by $\alpha_t$ and the brown arc $b$), $\partial D$, $\partial D'$, and $\partial D''$ are all white.}
    \label{Fi:Surgering_D}
\end{figure}

Surger ($\partial$-compress) $D$ along $X'$, as shown in Figure \ref{Fi:Surgering_D}.  This $\partial$-compression splits $D$ into two disks $D'$ and $D''$, each of which intersects $V$ fewer times than $D$ does, while $|\partial D'\cap L|+|\partial D''\cap L|=k$. The lexicographical minimality of $|\partial D\cap L|$ and $|D\cap V|$ implies that $\partial D_1$ and $\partial D_2$ are both contractible in $\partial(S\cut F)$, which contradicts $\partial D$ being homotopically essential in $\partial(S\cut F)$.  This contradiction completes the proof that, indeed, $\partial D_0,\dots,\partial D_\ell$ each contain at least two markers.
Since there are $k+2\ell$ markers in total, it follows that $\partial D_0,\dots,\partial D_\ell$ each contain at most $k$ markers. 

Finally, suppose that $k< n$. Then $D_0,\dots,D_\ell$ each have fewer than $n$ markers. Hence, each of $\partial E_0,\dots,\partial E_\ell$ is contractible in either $F_0$ or $F_1$, and thus in $F$. Ergo, for each $i=0,\dots,\ell$ there exists a smooth map $g_i$ from a disk $D^2_i$ to $F$ which restricts to an immersion $g_i:\partial D^2_i\to \partial E_i$. 

Glue the disks $D^2_0,\dots,D^2_\ell$ together into a single disk $D^*$ as follows.  For each $t=1,\dots,\ell$, there are exactly two indices $i,j$ (depending on $t$) for which $\beta_t\subset g_i(\partial D^2_i)$ and $\beta_t\subset g_j(\partial D^2_j)$.  Note that $g_i^{-1}(\beta_t)$ consists of an arc in $\partial D^2_i$ and possibly some isolated points, and likewise for $g_j^{-1}(\beta_t)$.  Glue $D^2_i$ and $D^2_j$ along these two arcs by an appropriate restriction of $g_j^{-1}\circ g_i$. Do this for each $t=1,\dots,\ell$ to obtain $D^*$.

Gluing the maps $g_0,\dots,g_\ell$ together in the same manner yields a continuous map $g:D^*\to F$ that satisfies $g(\partial D^*)=\partial D$.
Ergo, $D$ is $\partial$-contractible in $F$. Yet, this contradicts the assumption that $D\in\Caps_e(F)$. The conclusion of the theorem follows: $\ess(F)=|\partial D\cap L|=k\geq n$.
\end{proof}

    
We should note the close similarity of the preceding proof (until the last three paragraphs) and the proof of \cite[Theorem 7.3]{murasugi4D}, which asserts that any Murasugi sum of $\pi_1$-essential spanning solids in $S^4$ is $\pi_1$-essential (see \cite{murasugi4D} for the relevant definitions). That proof was adapted from the (overly complicated) proof of Theorem \ref{T:EndEss} in an earlier draft of this paper, and the proof that now appears above was adapted from the improved proof in \cite{murasugi4D}.

\section{Twisted plumbing, and the essence of an alternating checkerboard surface}\label{S:Twisted}

We have seen that the geometric and algebraic notions of essence sometimes behave differently---in particular, the latter behaves more naturally under Murasugi sum than the former, at least in general.  How aberrant are the cases where these notions differ?  We begin to address this question by focusing on a well-behaved class of surfaces: checkerboard surfaces from reduced alternating link diagrams. 

 \begin{notation}
Throughout \textsection\ref{S:Twisted}, $P\subset S^2$ is a reduced alternating diagram of a link $L$, $B$ is its positive-definite checkerboard surface, $W$ is its negative-definite checkerboard surface, $\Gamma$ is the Tait graph of $B$, and $G\in\Z^{\beta_1(B)\times\beta_1(B)}$ is a Goeritz matrix for $B$. 
\end{notation}

Do we always have $\ess(B)= \ess_g(B)$ for such a surface $B$?  Of course, we always have $\ess(B)\leq \ess_g(B)$, but how can we prove the reverse inequality in this setting?

Intuitively, one might predict that the essence---geometric or algebraic---of any such surface $B$ should equal the length of the shortest cycle in its Tait graph.  One can also express this intuition in terms of extreme framings of curves in these surfaces.  Shortly, we will validate this intuition in the case of geometric essence.  Verifying it in the algebraic setting, however, will require more work---this is where we will introduce a ``twisted'' generalization of Murasugi sum.  It will, of course, yield the corollary that $\ess(B)= \ess_g(B)$.

\begin{definition}
Write $\ess(\Gamma)$ for the length of the shortest cycle in $\Gamma$ and 
\[\ess(G)=\min\{\x^TG\x:~\text{nonzero }\x\in\Z^{\beta_1(B)}\}.\]
\end{definition}

\begin{prop}\label{P:ess_g(B)}
If $B$ is the positive-definite checkerboard surface for a reduced-alternating diagram $P$ of a link $L\subset S^3$, $\Gamma$ is a Tait graph of $B$, and $G$ is a Goeritz matrix for $P$, then
\[\ess(G)=\ess(\Gamma)=\ess_g(B).\]  
\end{prop}

An analogous statement holds for the negative-definite surface $W$ from $P$: the only difference is that the framings in the definition of $\ess(G)$ need to be negated.  Thus, the proposition states that the geometric essence of any checkerboard surface from any reduced alternating diagram equals the length of the shortest cycle in its Tait graph.  Similar comments apply throughout this section.

\begin{proof}[Proof of Proposition \ref{P:ess_g(B)}]
The equality $\ess(G)=\ess(\Gamma)$ follows from Theorem 5.1 of \cite{greene} (which Greene describes as a reframing of Theorem 1 of \cite{gordlith}), which states that the deformation retraction from $B$ to $\Gamma$ induces an isometry between $H_1(B)$ (under the Gordon-Litherland pairing, i.e. under $G$) and the flow lattice on $\Gamma$.  It will thus suffice to prove that $\ess_g(B)\leq \ess(\Gamma)$ and $\ess(G)\leq \ess_g(B)$.

Choose a cycle in $\Gamma$ of minimal edge-length $n$.  This yields a circle $\gamma\subset B$ that runs exactly once through exactly $n$ vertical arcs. Pushing $\gamma$ up through each vertical arc to the overpass gives a curve which bounds a geometric cap and intersects $L$ in $n$ points. Hence, $\ess_g(B)\leq \ess(\Gamma)$.

\begin{figure}[h]
\centering
\labellist
\pinlabel {$\longrightarrow$} at 118 45
\endlabellist
\includegraphics[width=.45\textwidth]{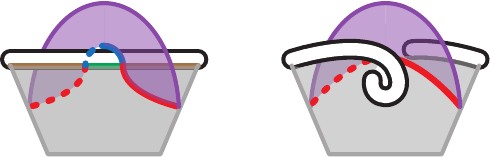}
\hfill
\labellist
\pinlabel {$\longrightarrow$} at 117 45
\endlabellist
\includegraphics[width=.45\textwidth]{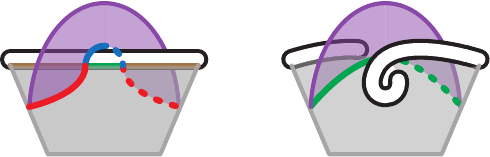}
\caption{Attaching a M\"obius band to $B$ near $D\cap L$ turns $D$ into a compressing disk}
\label{Fi:Cap_to_Compression}
\end{figure}

Consider a $\partial$-essential geometric cap $D$ for $B$ with $|\partial D\cap L|=\ess_g(B)=n$.  By assumption, $\partial D$ is essential in $B$; it is also unknotted, and so definiteness implies that $\partial D$ is not nullhomologous in $B$. This also holds for the circle $\gamma$ obtained by perturb $\partial D$ slightly into the interior of $B$.  We claim that $\gamma$ has framing at most $n/2$.  Indeed, if, instead of perturbing, we had replaced the disk-neighborhood in $B$ of each point of $\partial D\cap L$ with an appropriate M\"obius band as shown in Figure \ref{Fi:Cap_to_Compression}, then after these $n$ replacements, $\partial D$ would have been 0-framed in $B$; yet, each of these replacements should have changed the framing of $B$ near $\partial D$ by $\pm\frac{1}{2}$. Since the framing of $\gamma$ in $B$ is at most $n/2$, taking $\x\in \mathbb{Z}^{\beta_1(B)}$ to be a coordinate vector for $[\gamma]\in H_1(B)$ gives $\ess(G)\leq \x^TG\x\leq n=\ess_g(B)$.
\end{proof}

\subsection{Twisted Murasugi sum}\label{S:PPD}

The inquiry initiated above---is $\ess(B)=\ess_g(B)$?---now motivates a ``twisted'' generalization of Murasugi sum.  After providing a general definition, we immediately narrow our focus to instances that interact nicely (from a diagrammatic point of view) with the class of surfaces $B$ from above.

If a sphere $Q$ decomposes a spanning surface $F$ as a Murasugi sum $F=F_0*F_1$ along a plumbing disk $U=F\cap Q$, then the thickened sphere $\nu Q$ also intersects $F$ in a cobordism disk which is ``horizontal'' relative to the fibers of projection $\nu Q\to Q$.  The following notion generalizes Murasugi sum by adopting this perspective, but allowing cobordism disks that ``twist.''  

\begin{definition}\label{D:PP}
Let $F\subset S^3$ be a spanning surface and $Q\subset S^3$ a 2-sphere with regular neighborhood $\nu Q$.  View $\nu Q$ as $Q\times[-1,1]$ with $Q=Q\times\{0\}$, write $Q_\pm=Q\times\{\pm1\}$, and let $\pi_Q:\nu Q\to Q$ denote projection. Let $H_\pm$ denote the two components of $S^3\cut\nu Q$, and write $B_\pm=H_\pm\cup\nu Q$, $F_\pm=F\cap B_\pm$, and $U=F\cap\nu Q$.   Observe that $U$ is a 2-dimensional cobordism with boundary between the 1-manifolds with boundary $F\cap Q_+$ and $F\cap Q_-$.  

Isotope $F$ near $U$ so that, for each $x\in Q$, either $F\pitchfork \pi^{-1}(x)$ or $F\cap \pi^{-1}(x)$ is a properly embedded arc in $F$, doing this so as to minimize the number of points $x_1,\hdots,x_m\in Q$ over which the fiber of $\pi_Q$ is not transverse to $U$.  Write $x=\{x_1,\hdots,x_m\}$. Fixing $F\cap \pi_Q^{-1}(x)$, isotope $F$ near $U$ so that  $U\cap Q$  is a pinched disk with pinch points $x$.

Write $F=F_+\wt{*} F_-$, and call this a {\bf twisted Murasugi sum} along $U$ with {\bf defect} $m$ and {\bf complexity} $|\partial U\cap L|$. Call $U$ a {\bf twisted plumbing disk} for $F$, and call $Q\cut U$ an associated  (twisted plumbing) {\bf cap}.
There is an associated tree $T_U$ whose vertices and edges respectively correspond to the components of $U\cut x$ and the points of $x$. Call the edge-length-diameter of $T_U$ the {\bf diameter} of the twisted plumbing.
\end{definition}

We observe two basic properties of such a twisted Murasugi sum.

\begin{rem}\label{R:DM}
The diameter $d$ of a twisted Murasugi sum never exceeds its defect $m$: $d\leq m$.
\end{rem}

\begin{rem}\label{R:Pinch}
Given a twisted Murasugi sum along a twisted plumbing disk $U$ with pinch points $x$, consider a component $U_0$ of $U\cut x$. 
The number of pinch points on $\partial U_0$ must have the same parity as $|\partial U_0\cap L|$. Therefore, $|\partial U\cap L|$ is always even.  
\end{rem}

 We now narrow our focus to twisted Murasugi sums that play well with reduced alternating link diagrams and their checkerboard surfaces. (Note that a twisted Murasugi sum is a genuine Murasugi sum if and only if its defect is zero.)

\begin{figure}[h]
\begin{center}
\labellist
\small\hair 4pt
\pinlabel{$r-1$ crossings} at 270 50
\endlabellist
\includegraphics[height=1.75in]{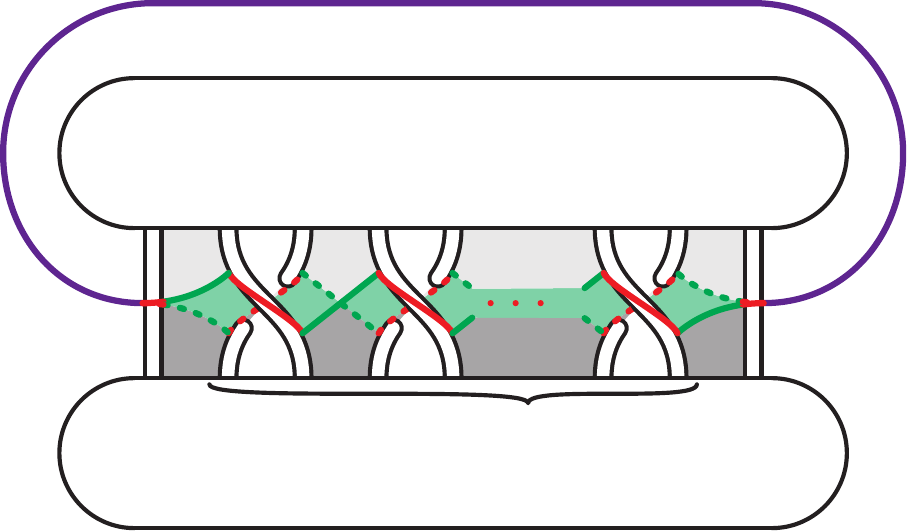}
\caption{A cap $D$ for a diagrammatic twisted plumbing decomposition of complexity $|\partial D\cap L|=2r$}
\label{Fi:DiagTwist}
\end{center}
\end{figure}

\begin{figure}[h]
\begin{center}
\labellist
\small \hair 4pt
\pinlabel {$=$} at 330 140
\pinlabel {$*$} at 500 180
\endlabellist
\includegraphics[height=2in]{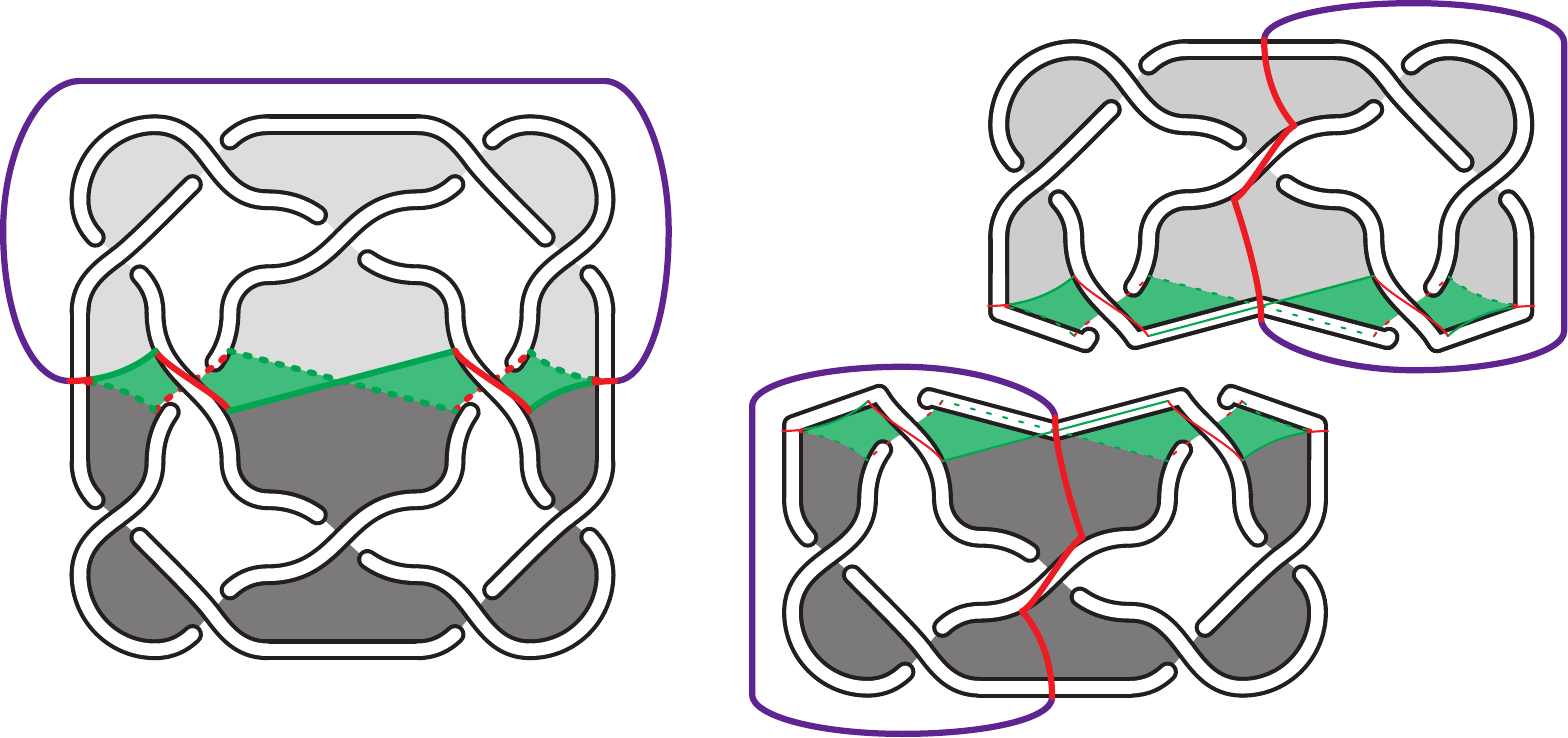}
\caption{A twisted deplumbing of a surface $F=F_0\widetilde{*}F_1$ with $\ess_c(F)=4=\ess_c(F_0)=\ess_c(F_1)$.}
\label{Fi:TwistHierarchy}
\end{center}
\end{figure}

\begin{definition}\label{D:DiagTwist}
Let $B$ and $W$ be the positive- and negative-definite checkerboard surfaces from a reduced alternating link diagram $P$, and let $D$ be an acceptable $\partial$-contractible cap for $B$ with $|D\cap W|=1$ and $|\partial D\cap L|=2r$.  Then $D$ appears as in Figure \ref{Fi:DiagTwist} and is called a {\bf diagrammatic plumbing cap} in case $r=1,2$ and a {\bf diagrammatic twisted plumbing cap} in general.  It defines a decomposition of $B=B_0\cup B_1$ into two subsurfaces where $B_0\cap B_1=U$.  Figure \ref{Fi:TwistHierarchy} shows an example.
\end{definition}

The $\partial$-contractible essence $\ess_c(B)$ (recall Definition \ref{D:Essence}) will play a more prominent role going forward. In the following observation, note that if $\beta_1(B)=1$ then $\ess_c(B)=\infty$:

\begin{prop}
Let $B$ be a checkerboard surface from a reduced alternating diagram.  If $\beta_1(B)\geq 2$, then $\ess_c(B)\leq 2(\ess(B)-1)$.  
\end{prop}

That is, $r+1=\frac{\ess_c(B)}{2}+1\leq \ess(B)\leq\ess_g(B)=n$.

\begin{figure}[h]
\begin{center}
\includegraphics[width=\textwidth]{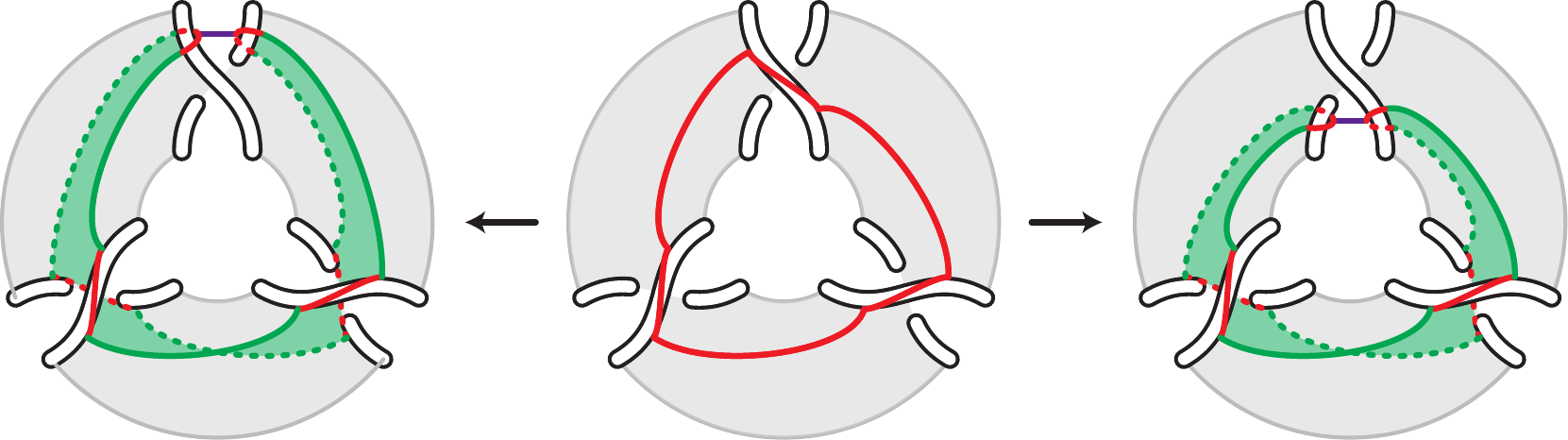}
\caption{A cycle of length $r$ in $\Gamma$ gives a twisted plumbing cap $D$ for $B$ of complexity $|\partial D\cap L|=2r$}
\label{Fi:Triangle}
\end{center}
\end{figure}

\begin{proof}
By Theorem \ref{T:CBEss}, the Tait graph $\Gamma$ from $B$ has a cycle of length $n=\ess(B)$, so there is an essential circle $\gamma\subset B$ that runs through exactly $n$ distinct crossings.  Then $B$ has a diagrammatic twisted deplumbing near $\gamma$ as shown in the case $n=3$ in Figure \ref{Fi:Triangle}. In the figure, either of the caps indicated left and right might be fake, but they cannot both be fake because $\beta_1(B)\geq 2$.
\end{proof}

\begin{prop}\label{P:EssC}
The $\partial$-contractible essence of $B$ equals the smallest complexity for which $B$ admits a \emph{diagrammatic} twisted plumbing cap,
\begin{equation}\label{E:EssC}
\text{ess}_c(B)=\min_{\text{diagrammatic }D\in\Caps(F)\setminus\Caps_e(F)}|\partial D\cap L|.
\end{equation}
\end{prop}

\begin{proof}
Every diagrammatic twisted plumbing cap is a $\partial$-contractible cap, so $\leq$ holds in (\ref{E:EssC}).  For the reverse inequality, consider a $\partial$-contractible cap $D$ for $B$ with $|\partial D\cap L|=\ess_c(B)=2r$.  Because $D$ is $\partial$-contractible, it must intersect $W$, so there are at least two outermost disks of $D\cut W$.  

Consider an outermost disk $D_0$ of $D\cut W$, and write $\alpha$ for the arc of $D\cap W$ in $\partial D_0$ and $\beta$ for the arc where $\partial D_0$ coincides with $\partial D$.  If the interior of $\alpha$ contains fewer than $r$ (entire) components of $\partial D\cap L$, then one can construct a diagrammatic twisted plumbing cap of complexity at most $2r$ as follows:  perturb $\alpha$ and $\beta$ so their endpoints lie on $L$, double $D_0$ across $S^2$, and modify near crossings as shown in Figure \ref{Fi:Cutting_path}.  In this case, we are done.

\begin{figure}[h]
\begin{center}
\labellist
\pinlabel {$\longrightarrow$} at 530 100
\endlabellist
\includegraphics[width=.475\textwidth]{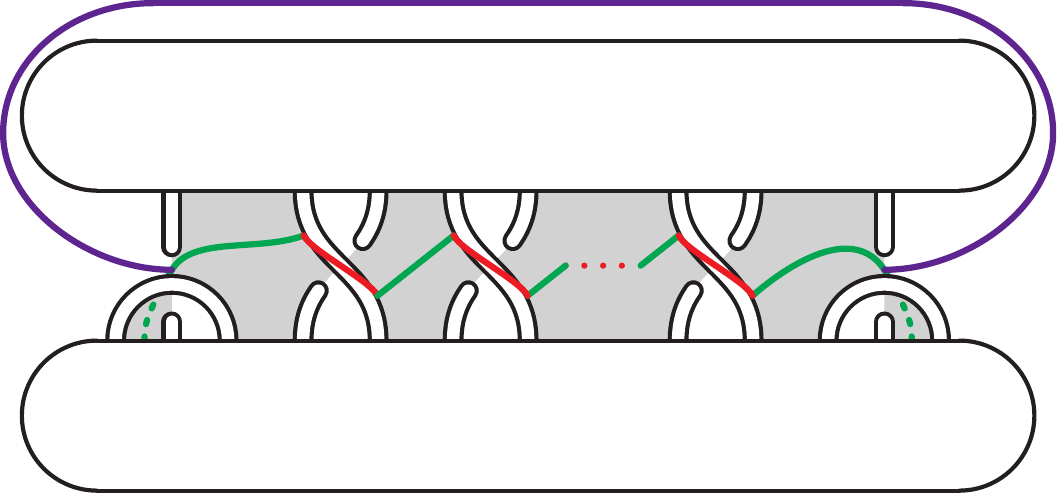}
\hfill
\includegraphics[width=.475\textwidth]{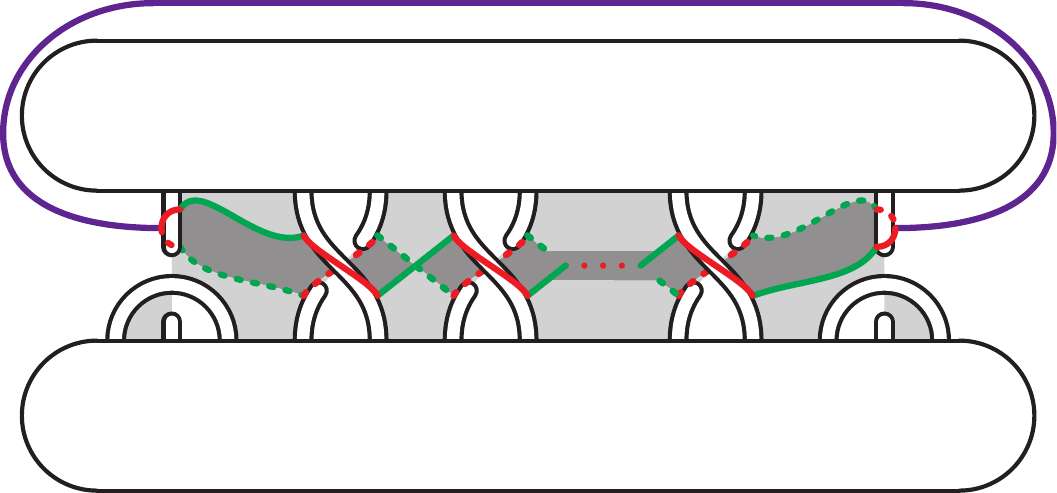}
\caption{Reading a minimal complexity diagrammatic plumbing cap off of a Tait graph}
\label{Fi:Cutting_path}
\end{center}
\end{figure}

The only remaining possibility is that $D\cap W=\alpha$, with both endpoints of $\alpha$ lying in $\mathring{B}$, rather than on $L$, but this is impossible because $\partial D$ must intersect each crossing band in $B$ an even number of times (since $D$ is $\partial$-contractible). 
\end{proof}

\begin{rem}\label{R:Read_off}
Recall from Remark \ref{R:Pinch} that quantity on either side of \eqref{E:EssC} is always even; we will often denote it $2r$.  One can read it off of the Tait graph $\Gamma$ for $B$, as $r-1$ is the edge-length is the shortest path $\rho$ in $\Gamma$ for which $\Gamma\setminus \rho$ is disconnected. (There is a straightforward correspondence between such separating paths in $\Gamma$ and the diagrammatic twisted deplumbings of $B$.)
\end{rem}

\begin{figure}[h]
\begin{center}
\includegraphics[width=.5\textwidth]{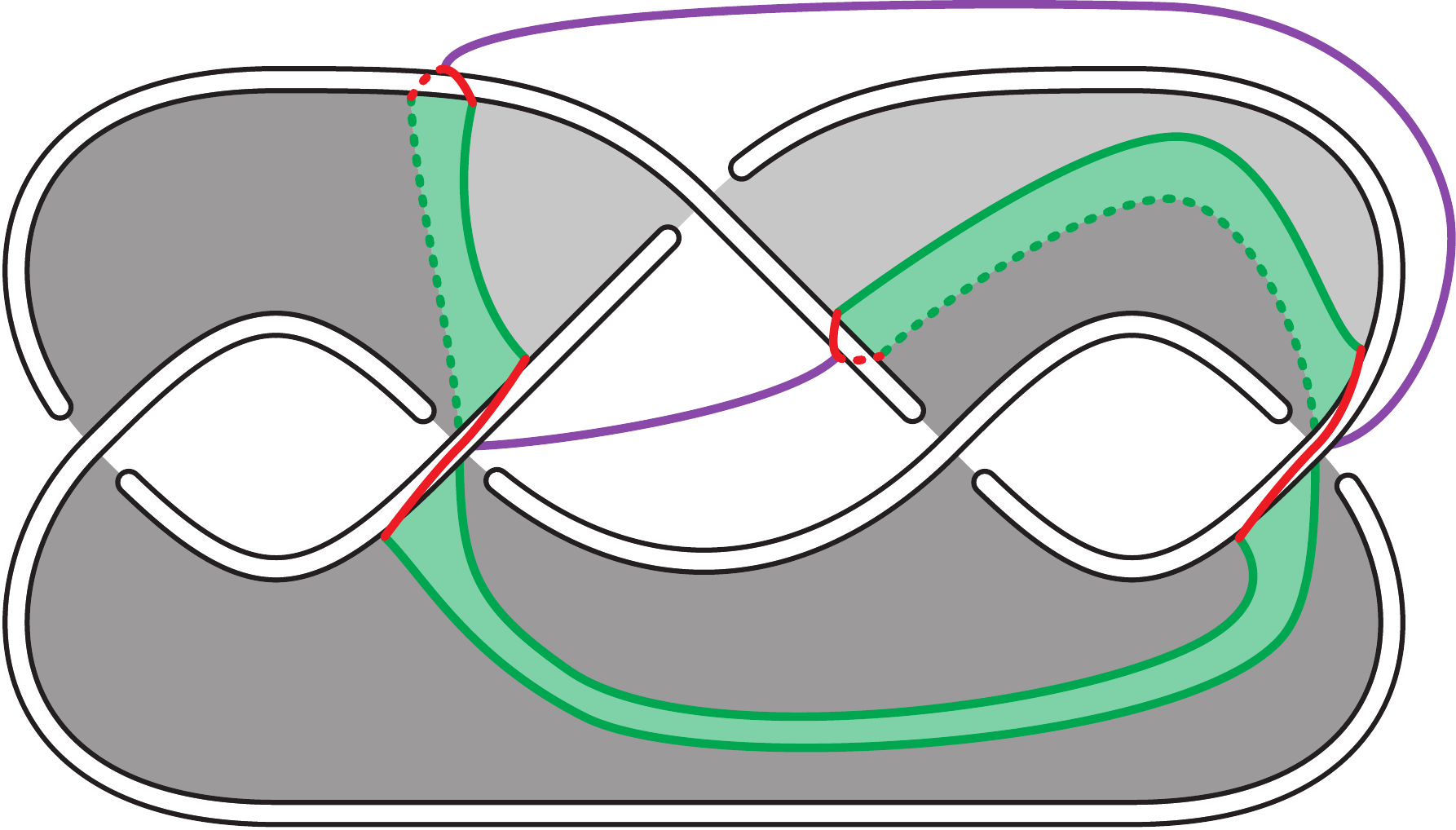}
\caption{A minimal-complexity plumbing cap which is not isotopic to the type of ``diagrammatic'' cap shown in Figure \ref{Fi:DiagTwist}}
\label{Fi:Whitehead}
\end{center}
\end{figure}

\begin{rem}\label{R:Non_diagrammatic}
It is not true, however, that \emph{all} minimal-complexity twisted deplumbings of a given surface $B$ are isotopic to diagrammatic ones---not even in the untwisted case.  Figure \ref{Fi:Whitehead} shows a counterexample. 
\end{rem}

Still, if one restricts to {\it diagrammatic} deplumbings of minimal complexity (which is possible because of Proposition \ref{P:EssC}), one can find a {well-behaved sequence} of twisted Murasugi sum decompositions of $B$:

\begin{definition}\label{D:Hierarchy}
Given a checkerboard surface $B$ from a reduced alternating link diagram, decompose $B=B_0\wt{*}B_1$ under diagrammatic twisted Murasugi sum of complexity $\ess_c(B)$.  Then decompose each $B_i=B_{(i,0)}\wt{*}B_{(i,1)}$ under diagrammatic twisted plumbing of minimal complexity $\ess_c(B_i)$.  Continue in this way.  

If one wishes to consider only decompositions with complexity less than some fixed amount $2r$, then one continues until each remaining surface $B_{\x}$ has $\ess_c(B_\x)\geq 2r$; otherwise, one continues until each $\beta_1(B_{\x})=1$ and thus $\ess_c(B_\x)=\infty$. 

Call this a {\bf hierarchical twisted deplumbing} of $B$.  One can recover $B$ from the resulting annuli and M\"obius bands $B_{\x}$ by gluing them together via twisted Murasugi sums in the opposite order.
\end{definition}

It is worth noting here that if $B=B_0\wt{*}B_1$ is a decomposition of minimal complexity, then it does not necessarily follow that both $\ess_c(B_i)\geq\ess_c(B)$.  Figure \ref{Fi:TwistedEssence} shows a counterexample.  Thus, the complexity of the twisted deplumbings may decrease during the sequence. On the other hand, it is not too hard to show that the set of factors that ultimately results from this decomposition is completely determined by $B$.  Because we will not need this fact, we leave the details to the reader.

\begin{prop}\label{P:(1)(2)}
Let $B$ be a checkerboard surface from a reduced alternating link diagram.  If $\ess_g(B)=n$ and $\ess_c(B)\geq 2(r+1)$ 
for $r=\lfloor\frac{n+1}{2}\rfloor$, then $\ess(B)=n$.
\end{prop}

\begin{proof}
By hypothesis, $n\geq 2$, and if $n=2$, the claim is just that $B$ is $\pi_1$-essential, which is well-known \cite{ozawa11}.  
Choose a $\partial$-essential cap $D$ for $B$ with $|\partial D\cap L|=\ess(B)$.   If $D\cap W=\varnothing$, then $D$ is also a geometric cap 
and we are done.  Assume otherwise.  Then each outermost disk of $D\cut W$ intersects $L$ at least $r$ times, not counting endpoints of arcs of $D\cap W$, because $\ess_c(B)\geq 2r$ (here we are using Remark \ref{R:Read_off} and Figure \ref{Fi:Cutting_path}). Thus, $\ess(B)\geq 2r\geq n$.  This completes the proof because $\ess(B)\leq \ess_g(B)$ holds in general.
\end{proof}

\subsection{Twisted plumbing respects sufficiently large essence}

The proof of the theorem below will use the contrapositive of Proposition \ref{P:FakeString}, which we state here for convenience.

\begin{cor}\label{C:FakeString}
If $D$ is a cap or fake cap for $F$ with $|\partial D\cap L|=2r$, and if $\partial D$ contains a string of $r$ consecutive essential arcs in $F$, then $D$ cannot be fake.
\end{cor}

\begin{theorem}\label{T:TwistedEss}
Assume $F=F_0\wt{*}F_1$ is a twisted Murasugi sum along a twisted plumbing disk $U$ and an associated cap $V$, where $|\partial V\cap L|=2r=\ess_c(F)$. If $\min_{i=0,1}\ess(F_i)=n\geq 2r$, then $\ess(F)\geq n$.
\end{theorem}

Note that, unlike in Theorem \ref{T:PlumbEss}, the complexity of the twisted plumbing is bounded here in terms of the essences of the surfaces: $n\geq 2r$.  This will allow us to use an outermost disk argument, rather than a more complicated counting argument like the one in the proof of Theorem \ref{T:PlumbEss}. Consequently, any concerns that arise about possibly fake caps will be easier to resolve. (Also see Remark \ref{R:TwistedSharp}.) Also, we remind the reader that caps are assumed to be non-fake unless stated otherwise.

\begin{proof}
The assumption that $|\partial V\cap L|=2r=\ess_c(F)$ implies that $V$ is acceptable and in particular that each arc of $\partial V\cap F$ is essential in $F$.
Assume for contradiction that $D$ is a $\partial$-essential cap for $F$ such that $|\partial D\cap L|=n'<n$, and, subject to this and the generic assumption that $\partial D\cap L\cap\partial V=\varnothing$, that $(|D\cap V|,|\partial D\cap L|)$ has also been minimized lexicographically. Then $D$ must intersect $V$ and $|\partial V\cap L|\leq 2r$, so there is an outermost disk $V_0\subset V\cut D$ whose boundary consists of an arc $\alpha$ of $D\cap V$ and an arc $\alpha'$ of $\partial V\cut\partial D$ with $|\alpha'\cap L|=r'\leq r$.  Surger $D$ along $V_0$ to get two caps (possibly fake) $D'$ and $D''$ for $F$. 

Assume without loss of generality that $|\partial D'\cap L|\geq |\partial D''\cap L|$. Then $|\partial D'\cap L|>|\partial D\cap L|$, or else minimality would imply that both $D'$ and $D''$ are $\partial$-contractible, which would apply contrary to assumption that $D$ is $\partial$-contractible. 
Also, $|D''\cap V|<|D\cap V|$, so minimality implies that $D''$ is not $\partial$-essential.  
Further,
\[|\partial D''\cap L|=|\partial D\cap L|-|\partial D'\cap L|+2r'<2r'\leq 2r=\ess_c(F),\]
so $D''$ cannot be a $\partial$-contractible cap either. Instead, it must be fake. Thus, Observation \ref{O:FakeDegree} implies that $|\partial D''\cap L|$ is even, hence that $|\partial D''\cap L|\leq 2r'-2$.  This contradicts Corollary \ref{C:FakeString}, because $\partial D''$ contains the string of $r'-1$ arcs of $\partial V\cut L$ that comprise $\alpha'\cut L$, and all of these arcs are essential in $F$.
\end{proof}

\begin{figure}[h]
\begin{center}
\labellist
\small \hair 4pt
\pinlabel {$=$} at 330 140
\pinlabel {$*$} at 480 160
\endlabellist
\includegraphics[height=2in]{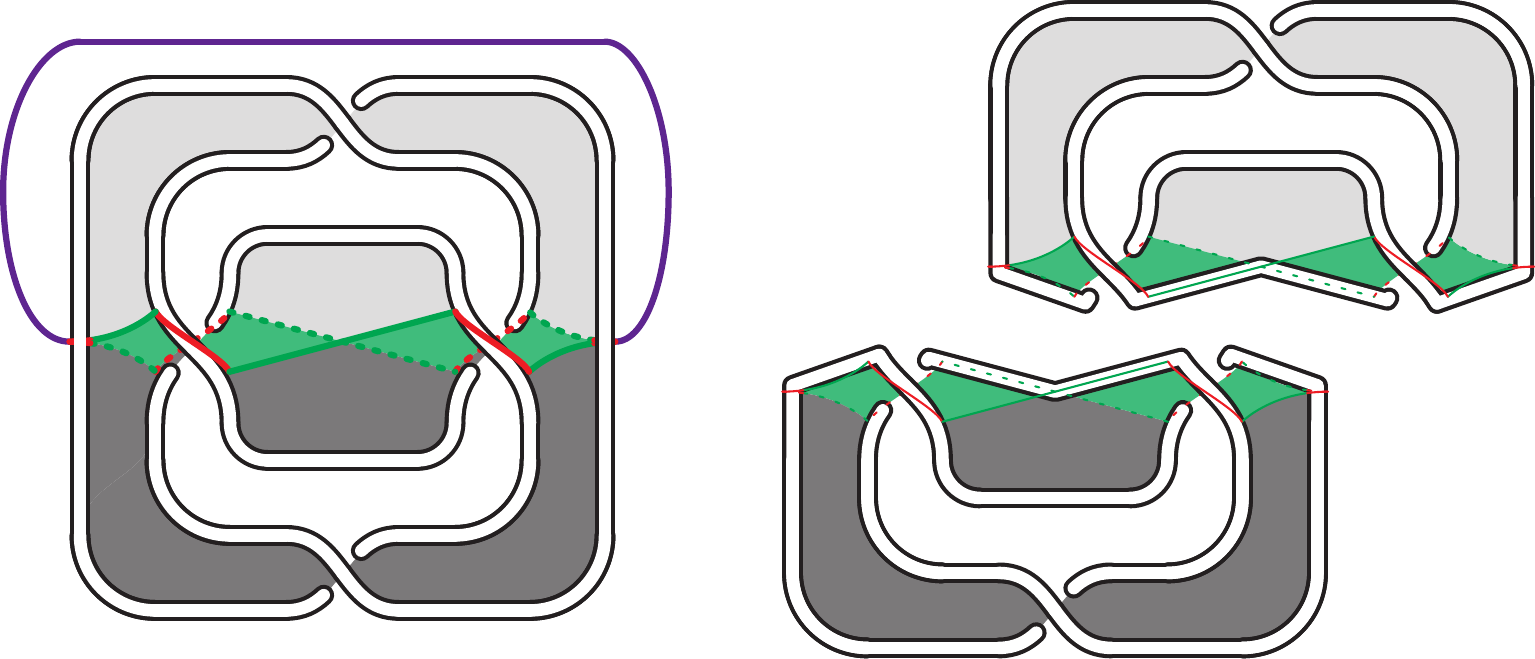}
\caption{A twisted deplumbing $F=F_0\widetilde{*}F_1$ with $\ess(F_0)=3=\ess(F_1)$ and $\ess(F)=2$}
\label{Fi:TwistedEssence}
\end{center}
\end{figure}

\begin{rem}\label{R:TwistedSharp}
We chose the bound $n\geq 2r$ in Theorem \ref{T:TwistedEss} because this simplifies the proof of this theorem and it will be sufficient when we apply this result in the proof of Theorem \ref{T:CBEss}.  It seems unlikely, however, that this bound is sharp.  Some bound, however, is necessary.  For example, if $F_0$ and $F_1$ are both M\"obius bands spanning $(2,n)$ torus knots, each drawn as a checkerboard surface from a reduced alternating diagram of $n$ crossings, and if $F=F_0\wt{*}F_1$ is a diagrammatic twisted Murasugi sum of complexity $n+3$, then $\ess(F)=n-1$. See Figure \ref{Fi:TwistedEssence}. 
\end{rem}

This raises the following natural question:

\begin{question}
Does Theorem \ref{T:TwistedEss} hold under the weaker assumption that $n\geq 2r-2$?\end{question}

\subsection{Essence of a state surface from an alternating diagram}%

\begin{theorem}\label{T:CBEss}
If $B$ is the positive-definite checkerboard surface from a reduced alternating link diagram $P$, then $\ess(B)=\ess_g(B)$.
\end{theorem}


\begin{proof}
Following Propositions \ref{P:ess_g(B)} and \ref{P:(1)(2)}, write $n=\ess_g(B)=\ess(\Gamma)$ and $r=\lfloor \frac{n+1}{2}\rfloor$, so that 
$2r\in \{n,n+1\}$.
Hierarchically decompose  $B=B_1\wt{*}\cdots \wt{*}B_t$ under diagrammatic twisted deplumbings of complexities at most $2r$, as described in Definition \ref{D:Hierarchy} (if $\ess_c(B)> 2r$, this leaves $B$ unchanged). 

Then each $B_i$ is the positive-definite checkerboard surface from a reduced alternating link diagram and its Tait graph $\Gamma_i$ is a subgraph of $\Gamma$. Hence, writing $\ess(\Gamma_i)=n_i$ for each $i$, we have $n_i\geq n$. 
Also, each $\ess_c(B_i)> 2r$ by assumption and is even by Remark \ref{R:Pinch}, so Proposition \ref{P:(1)(2)} further implies that each $\ess(B_i)\geq n$. 

Finally, since each twisted Murasugi sum in the reconstruction of $B=B_1\wt{*}\cdots \wt{*}B_t$ has complexity at most $n=2r$ and yields a surface whose $\partial$-contractible essence equals the complexity of the last plumbing performed, Theorem \ref{T:TwistedEss} implies inductively that $\ess(B)\geq n$. This completes the proof because $\ess(B)\leq \ess_g(B)$ holds in general.
\end{proof}

\begin{cor}\label{C:StateEss}
Let $x$ be a 
 homogeneously adequate state, $F_x$ its state surface (with any layering of the state disks), and $\Gamma_x$ its state graph. Write $\ess(\Gamma_x)$ for the length of the shortest cycle in $\Gamma_x$. Then $\ess_g(F_x)\geq \ess(F_x)\geq \ess(\Gamma_x)$. Moreover, $x$ has a state surface $F'_x$ (possibly with a different layering of the state disks) satisfying $\ess_g(F'_x)=\ess(F'_x)=\ess(\Gamma_x)$.
\end{cor}

\begin{proof}
Remark \ref{R:CBPlumb} explains why the first claim follows from Theorems \ref{T:CBEss} and \ref{T:PlumbEss}. For the second, choose a shortest cycle in $\Gamma_x$.  It runs through a single checkboard plumbing factor of $F_x$.  Place all state disks in that factor surface above the projection sphere $S^2$, and place the rest of the state disks below $S^2$.  The needed cap now lies entirely above $S^2$.
\end{proof}

\section{Plumbing end-essential spanning surfaces in arbitrary 3-manifolds}\label{S:3Mfld}

Finally, we extend everything from $S^3$ to the setting of an arbitrary 3-manifold.

\begin{convention}
Throughout \textsection\ref{S:3Mfld},
$F$ is a spanning surface for a link $L$ in a compact 3-manifold $M$ (not necessarily closed or orientable). View $F$ as a compact surface in $M$ with $\partial F=L$. Also, 
use the quotient map $h_F:M\cut F\to M$, which reglues corresponding pairs of points from $\mathring{F}$ in $\partial (M\cut F)$, to write $\wt{L}=h_F^{-1}(L)$, $\wt{F}=h_F^{-1}(\mathring{F})$, and $\wt{\partial M}=h_F^{-1}(\partial M)$.  
\end{convention}

\begin{definition}\label{D:EndEss}
With the setup above, say that $F$ is (algebraically):
\begin{enumerate}
\item {\bf incompressible} if any circle
$\gamma\subset \wt{F}$ that bounds a disk $\wt{D}\subset M\cut F$  
also bounds a disk in $\wt{F}$. In that case, $D=h_F(\wt{D})$ is called a {\it fake compressing disk} for $F$; if $\gamma$ does not bound a disk in $\wt{F}$ then $D$ is a {\it compressing disk}.
\item {\bf end-incompressible} if 
any circle %
$\gamma\subset \wt{F}$ that is parallel through an annulus $\wt{A}$ in $M\cut F$ to $\wt{\partial M}$ bounds a disk in $\wt{F}$. In that case, $A=h_F(\wt{A})$ is called a {\it fake end-annulus} for $F$; if $\gamma$ does not bound a disk in $\wt{F}$ then $A$ is an {\it end-annulus} for $F$.
\item {\bf $\bs{\partial}$-incompressible} if, for any circle
$\gamma\subset \partial (M\cut F)$ with $|\gamma\cap\wt{L}|=1$ that bounds a disk $\wt{D}$ in $M\cut F$, the arc $\gamma\cut\wt{L}$ is parallel in $\partial (M\cut F)\cut\wt{L}$ to $\wt{L}$. If $\gamma\cut\wt{L}$ is not so parallel, then $h_F(\wt{D})$ is a {\it $\partial$-compressing disk} for $F$.
\item {\bf $\pi_1$-essential} if $F$ satisfies (a) and (c).
\item {\bf end-essential} if $F$ satisfies (a), (b) and (c).
\end{enumerate}
\end{definition}

This definition first appeared in \cite{endess}, where $L$ is a link in a thickened closed orientable surface  $M=\Sigma\times I$ and $P$ is a diagram of $L$ on $\Sigma$.  We will indicate this setting by writing pairs $(\Sigma,L)$ and $(\Sigma,P)$. Then we will sometimes also assume that $P$ is alternating on $\Sigma$, or that $P$ is {\bf cellular}, meaning that it cuts $\Sigma$ into disks.  

A crossing $c$ in $P$ is {\bf removably nugatory} if there is a disk $D\subset \Sigma$ such that $\partial D\cap P=\{c\}$. In that case, one can remove $c$ from $P$ by using a flype to move $c$ across the part of $P$ in $D$ and then using a Reidemeister 1 move to undo the resulting monogon. If $P$ has a removably nugatory crossing, then $B$ or $W$ is $\partial$-compressible.  The main result of \cite{endess} is the following strong converse of that fact:

\begin{theorem}[Theorem 1.1 of \cite{endess}]\label{T:EndEss}
 If $M=\Sigma\times I$ is a thickened orientable surface and $P$ is a cellular alternating link diagram on $\Sigma$ without removably nugatory crossings, then both checkerboard surfaces from $P$ are end-essential in $M$.
\end{theorem}

We will extend this theorem so that its conclusion holds for all adequate state surfaces from $P$---the state surface construction in $\Sigma\times I$ is the same as in the classical setting, with the extra requirement that each state circle must bound a disk in $\Sigma$, or else one could not cap it off with a disk in $\Sigma\times I$. To do this, we will extend Ozawa's theorem to plumbings of spanning surfaces in 3-manifolds other than $S^3$. First, we generalize 
Definition \ref{D:Plumb} as follows:

\begin{definition}\label{D:Plumb3Mfld}
Let $V\subset M$ be an embedded disk with $V\cap F=\partial V$ such that $\partial V$ bounds a disk $U\subset F$. Then the 2-sphere $Q=U\cup V$ separates $M$, hence decomposes $M$ as a (possibly trivial) connect sum, $M=M_0\# M_1$, where each $M_i$ is constructed by attaching a 3-handle to a component $N_i$ of $M\cut Q$. If neither subsurface $F_i=F\cap N_i$ is a disk, then $V$ is a {{\bf plumbing cap}} for $F$, and $F$ is obtained by {\bf plumbing} $F_0$ and $F_1$ along $U$, denoted $F_0*F_1=F$.  We may also call this a {\bf Murasugi sum}.
\end{definition}

The definitions of caps, acceptable caps, essence, and twisted plumbing (\ref{D:Caps}, \ref{D:Acceptable}, \ref{D:Essence}, and \ref{D:PP}) extend to the general 3-manifold setting without further comment.

\begin{theorem}\label{T:PlumbEndEss}
Suppose $M=M_0\#M_1$ is a (possibly trivial) connect sum of 3-manifolds and $F=F_0*F_1$ is a Murasugi sum of $\pi_1$-essential spanning surfaces $F_i\subset M_i$. Write $\min_{i=0,1}\text{ess}(F_i)=n$. Then $F$ is $\pi_1$-essential with $\ess(F)\geq n$. Moreover, if $M_0$ is a 3-sphere, $M_1$ is a thickened surface, and $F_1$ is end-essential, then $F$ is also end-essential.
\end{theorem}

\begin{proof}
The proof the first claim is the same as the proof of Theorem \ref{T:PlumbEss}, where $U$ is a plumbing disk for $F=F_0*F_1$, $V$ is an associated cap, and $D\in \Caps_e(F)$ is chosen so as lexicographically to minimize $k=|\partial D\cap\partial F|$, $\ell=|D\cap V|$, and $m=|\partial D\cap\partial U|$. The key observation is that each arc $\alpha_j$ of $D\cap V$ is once again parallel through each 3-handle $M_i\cut N_i$ to an arc $\beta_j\subset U$. 

The proof of the second claim is similar: if $F$ is not end-essential, choose an end-annulus $A$ for $F$ which lexicographically minimizes $(|A\cap V|,|\partial A\cap \partial U|)$.  Then $A$ must intersect $V$ in some nonempty collection of arcs $\alpha_1,\hdots, \alpha_\ell$, each with both endpoints on $F$, so $V$ cuts $A$ into $\ell+1$ pieces: $\ell$ disks and an annulus (the annulus contains the circle $\partial A\cap\partial M$).  Mark each endpoint of each arc $\alpha_i$ as before. (There is no need for the other type of marker, because $\partial A\cap L=\varnothing$.) With $2\ell$ markers and $\ell+1$ components of $A\cut V$, some component $A_0$ has fewer than two markers.  It extends as before to a (possibly fake) compressing disk, $\partial$-compressing disk, or end-annulus $Z$ for $F_0$ (or $F_1$).  There are no fake $\partial$-compressing disks, so $Z$ must be a fake compressing disk or a fake end-annulus.

Like before, carefully lift $\partial Z\cap F$ to a circle $\gamma\subset \partial (N_0\cut F_0)$, which bounds a disk $Y\subset \partial (N_0\cut F_0)$ because $Z$ is fake.  The lift of $\partial U$ in $\partial (N_0\cut F_0)$ cuts $Y$ into subdisks. Take an outermost subdisk, and consider the corresponding disk $Y_*$ of $(F_0\cup V)\cut\partial U$, whose boundary consists of an arc $\sigma\subset\partial U$ and an arc $\tau\subset (F_0\cup V)\cut\partial U$.  Figure \ref{Fi:Trivial_Disk_F} again shows the possibilities for $Y_*$, and again each contradicts minimality: the first two by a small isotopy of $\tau\subset\partial A$ past $\sigma\subset\partial U$, and the third by surgering $A$ along $Y_*$, after which minimality implies as before that $\partial A$ is contractible in $F$, and now, since $\partial A\cap L=\varnothing$, this implies contrary to assumption that $A$ is fake.
\end{proof}

Combining Theorems \ref{T:EndEss} and \ref{T:PlumbEndEss} yields:

\begin{theorem}\label{T:Endess}
If $P\subset \Sigma$ is a cellular alternating diagram without removably nugatory crossings, then every adequate state surface from $P$, under any layering of its state disks, is end-essential in $\Sigma\times I$.
\end{theorem}

\begin{proof}
Given a state surface $F_x$ and a non-innermost state circle $x_0$ of the underlying state $x$, take $U$ to be the disk that $x_0$ bounds in $F_x$, and take $V$ to be a disk with interior disjoint from $F_x$ that $x_0$ bounds on the opposite side of $\Sigma$. Then the 2-sphere $U\cup V$ bounds a ball $B_1\subset \Sigma\times I$, and $F_x\cap B_1$ is nontrivial because $x_0$ is non-innermost and $P$ has no removably nugatory crossings. Therefore, $F_x$ decomposes along the collection of all such spheres as a plumbing of a checkerboard surface $F$ for a cellular alternating diagram on $\Sigma$ and several checkerboard surfaces $F_i$ for cellular alternating diagrams on $S^2$.  Moreover, none of these diagrams have removably nugatory crossings, so, by Theorem \ref{T:EndEss}, all the $F_i$ are $\pi_1$-essential, and $F$ is end-essential.  Therefore, by Theorem \ref{T:PlumbEndEss}, $F_x$ is end-essential as well.
\end{proof}

The same argument proves more generally:

\begin{theorem}\label{T:EndEssH}
Every homogeneously adequate state surface in a thickened surface is end-essential.
\end{theorem}



\begin{thebibliography}{199}


\bibitem
{ak} C. Adams, T. Kindred, {\it A classification of spanning surfaces for alternating links}, Alg. Geom. Topol. 13 (2013), no. 5, 2967-3007.
%
%
%
%
%
\bibitem
{cklsv2} M. Cohen, T. Kindred, A. Lowrance, A. Shanahan, C. Van~Cott, {\it Surgery moves connecting all essential spanning surfaces for a 2-bridge link I}, Forthcoming.
%
\bibitem
{cklsv3} M. Cohen, T. Kindred, A. Lowrance, A. Shanahan, C. Van~Cott, {\it Surgery moves connecting all essential spanning surfaces for a 2-bridge link II}, Forthcoming.
%
%
%
%
%
\bibitem
{gab1} D. Gabai, {\it The Murasugi sum is a natural geometric operation},  Low-dimensional topology (San Francisco, Calif., 1981), 131-143, Contemp. Math., 20, Amer. Math. Soc., Providence, RI, 1983.
%
\bibitem
{gab2} D. Gabai, {\it The Murasugi sum is a natural geometric operation II},  Combinatorial methods in topology and algebraic geometry (Rochester, N.Y., 1982), 93-100, Contemp. Math., 44, Amer. Math. Soc., Providence, RI, 1985. 
%
%
%
%
%
%
\bibitem
{gordlith} C. McA. Gordon, R.A. Litherland, {\it On the signature of a link}, Invent. Math. 47 (1978), no. 1, 53-69.
%
%
\bibitem
{greene} J. Greene, {\it Alternating links and definite surfaces}, with an appendix by A. Juhasz, M Lackenby, Duke Math. J. 166 (2017), no. 11, 2133-2151.
%
%
\bibitem
{ht} A. Hatcher, W. Thurston, {\it Incompressible surfaces in 2-bridge knot complements}, Inv. Math. 79 (1985), 225-246.
%
%
%
%
%
%
\bibitem
{murasugi4D} H. Karimi, S. Kim, T. Kindred, M. Miller, P. Naylor, {\it Spanning solids and Murasugi sum in dimension 4}, \url{https://arxiv.org/abs/2607.10018}
%
%
%
%
%
%
\bibitem
{TkThesis} T. Kindred, {\it Checkerboard plumbings}, PhD thesis, University of Iowa (2018), Available at \url{https://
doi.org/10.17077/etd.yv1nlfcz}
%
\bibitem
{flyping} T. Kindred, {\it A geometric proof of the flyping theorem}, Adv. Math. 482 (2025) C, 110642, 53 pp.
%
%
\bibitem
{endess} T. Kindred, {\it End-essential spanning surfaces for links in thickened surfaces}, Topol. Appl. 361 (2025), 109187.
%
%
\bibitem
{natural} T. Kindred, {\it How natural of a geometric operation is Murasugi sum?}, Forthcoming.
%
%
%
%
%
%
%
%
%
%
%
%
%
\bibitem{mur63} K. Murasugi, {\it On a certain subgroup of the group of an alternating link}, Amer. J. Math. 85 (1963), 544-550.
%
%
%
\bibitem
{ozawa11} M. Ozawa, {\it Essential state surfaces for knots and links}, J. Aust. Math. Soc. 91 (2011), no. 3, 391-404. 
%
\bibitem{pardon11} J. Pardon, {\it On the distortion of knots on embedded surfaces}, Ann. of Math. (2) 174 (2011), no. 1, 637-646. 
%
%
%
%
%
%
%
%
%
%
%
%
%
\end{thebibliography}
\end{document}